\documentclass[11pt]{amsart}

\usepackage{amsfonts}
\usepackage{mathrsfs}
\usepackage{bbm}
\usepackage{amsmath}
\usepackage{amssymb}
\usepackage{amscd}
\usepackage{amsthm}

\makeatletter
\@namedef{subjclassname@2020}{\textup{2020} Mathematics Subject Classification}
\makeatother

\newtheorem{thm}{Theorem}[section]
\newtheorem{lemma}[thm]{Lemma}
\newtheorem{prop}[thm]{Proposition}
\newtheorem{cor}[thm]{Corollary}

\numberwithin{equation}{section}

\newcommand{\te}{\textstyle}

\newcommand{\ve}{\varepsilon}

\newcommand{\bc}{\mathbb{C}}
\newcommand{\bff}{\mathbb{F}}
\newcommand{\bh}{\mathbb{H}}
\newcommand{\bn}{\mathbb{N}}

\newcommand{\br}{\mathbb{R}}
\newcommand{\bz}{\mathbb{Z}}

\newcommand{\sa}{\mathscr{A}}

\newcommand{\sce}{\mathscr{E}}

\newcommand{\scg}{\mathscr{G}}
\newcommand{\ci}{\mathscr{I}}

\newcommand{\su}{\mathscr{U}}

\newcommand{\sd}{\mathrm{D}}

\newcommand{\cs}{\mathscr{S}}

\newcommand{\ca}{\mathcal{A}}
\newcommand{\cb}{\mathcal{B}}
\newcommand{\cd}{\mathcal{D}}
\newcommand{\cf}{\mathcal{F}}

\newcommand{\ch}{\mathcal{H}}

\newcommand{\cm}{\mathcal{M}}

\newcommand{\fg}{\mathfrak{g}}

\newcommand{\fv}{\mathfrak{V}}
\newcommand{\fw}{\mathfrak{W}}

\newcommand{\tv}{V^{\ast}}

\newcommand{\autg}{Aut\left(G \right)}
\newcommand{\autag}{Aut_a\left(G \right)}

\newcommand{\supp}{\mathrm{supp}}

\newcommand{\fvm}{\mathcal{F}(V, M)}

\newcommand{\ftscgt}{\tilde{\mathcal{F}}(\pi_{\mathscr{G}_2, W})}

\newcommand{\fvg}{\mathcal{F}\left(V, G \right)}
\newcommand{\fwg}{\mathcal{F}\left(W, G \right)}

\newcommand{\fve}{\mathcal{F}\left(V, \mathscr{E} \right)}

\newcommand{\fcve}{\mathcal{F}_c \left(V, \mathscr{E} \right)}

\newcommand{\fcvfg}{\mathcal{F}_c\left(V, \fg \right)}

\newcommand{\sg}{\mathcal{M}\left(S, G \right)}

\newcommand{\msg}{\mathcal{M}\left(S, G \right)}

\newcommand{\kp}{W^{k, p}}

\newcommand{\ka}{C^{k, \alpha}}

\newcommand{\fvrn}{\mathcal{F}\left(V, \mathbb{R}^n \right)}
\newcommand{\fvr}{\mathcal{F}\left(V, \mathbb{R} \right)}

\newcommand{\fvrp}{\mathcal{F}\left(V, \mathbb{R}^+ \right)}

\newcommand{\fczvrp}{\mathcal{F}_c^0\left(V, \mathbb{R}^+ \right)}

\newcommand{\ckovr}{C^{k_1}\left(V, \mathbb{R} \right)}
\newcommand{\ckovrp}{C^{k_1}\left(V, \mathbb{R}^+ \right)}
\newcommand{\ckvr}{C^{k}\left(V, \mathbb{R} \right)}

\newcommand{\fwr}{\mathcal{F}\left(W, \mathbb{R} \right)}
\newcommand{\tfwr}{\tilde{\mathcal{F}}\left(W, \mathbb{R} \right)}
\newcommand{\fwrn}{\mathcal{F}\left(W, \mathbb{R}^n \right)}

\newcommand{\tfwrp}{\tilde{\mathcal{F}}\left(W, \mathbb{R}^+ \right)}

\newcommand{\cktwr}{C^{k_2}\left(W, \mathbb{R} \right)}
\newcommand{\ckwr}{C^{k}\left(W, \mathbb{R} \right)}

\newcommand{\fcvr}{\mathcal{F}_c\left(V, \mathbb{R} \right)}

\newcommand{\fvf}{\mathcal{F}\left(V, \mathbb{F} \right)}

\newcommand{\fcvf}{\mathcal{F}_c\left(V, \mathbb{F} \right)}

\newcommand{\fcvg}{\mathcal{F}_{c}\left(V, G \right)}

\newcommand{\fcvscg}{\mathcal{F}_{c}\left(\pi_{\scg, V} \right)}

\newcommand{\fvscg}{\mathcal{F}\left(\pi_{\scg, V} \right)}
\newcommand{\fvscgo}{\mathcal{F}\left(\pi_{\scg_1, V} \right)}

\newcommand{\fcvx}{\mathcal{F}_{c}\left(\pi_{\xi, V} \right)}

\newcommand{\fczvg}{\mathcal{F}^0_{c}\left(V, G \right)}

\newcommand{\fczvgo}{\mathcal{F}^0_{c}\left(V, G_1 \right)}

\newcommand{\fczvscg}{\mathcal{F}^0_{c}\left(\pi_{\scg, V} \right)}

\newcommand{\fcvscgo}{\mathcal{F}_{c}\left(\pi_{\scg_1, V} \right)}
\newcommand{\fczvscgo}{\mathcal{F}^0_{c}\left(\pi_{\scg_1, V} \right)}

\newcommand{\tx}{\tilde{x}}

\newcommand{\dimv}{\dim_{\br} V}

\newcommand{\zero}{\mathbf{0}}
\newcommand{\one}{\mathbf{1}}

\newcommand{\fh}{\hat{f}}
\newcommand{\xh}{\hat{x}}

\newcommand{\ah}{\hat{a}}
\newcommand{\bhat}{\hat{b}}

\newcommand{\vb}{\bar{v}}

\newcommand{\ax}{x^{\ast}}

\newcommand{\xv}{\mathfrak{X}_c(V)}
\newcommand{\xw}{\mathfrak{X}_c(W)}
\newcommand{\citv}{\mathfrak{X}(V)}
\newcommand{\citw}{\mathfrak{X}(W)}

\newcommand{\foceo}{\cf_c(\pi_{\eta_1, V})}

\newcommand{\ftet}{\tilde{\cf}(\pi_{\eta_2, W})}

\newcommand{\zsao}{Z(\mathscr{A}_1)}

\newcommand{\fsao}{\cf(\pi_{\mathscr{A}_1, V})}
\newcommand{\fcsao}{\cf_c(\pi_{\mathscr{A}_1, V})}

\newcommand{\fsat}{\cf(\pi_{\mathscr{A}_2, W})}
\newcommand{\fcsat}{\cf_c(\pi_{\mathscr{A}_2, W})}

\newcommand{\ftsat}{\tilde{\cf}(\pi_{\mathscr{A}_2, W})}

\begin{document}

\title{Range decreasing group homomorphisms and weighted composition operators}

\author{Ning Zhang}

\address{School of Mathematical Sciences \\ Ocean University of China \\Qingdao, 266100, P. R. China}

\email{nzhang@ouc.edu.cn}



\thanks{The author is grateful to L. Lempert for his
very helpful suggestions.
This research was partially supported by the Natural Science Foundation of Shandong Province of China
grant ZR2023MA035.}

\begin{abstract}
We present necessary and sufficient conditions for a group homomorphism
between spaces of smooth sections of Lie group bundles
to be a weighted composition operator.
These results provide new insights into a wide range of problems related to weighted composition operators.
Specifically, we prove that the algebraic structure
of the space of smooth sections of an algebra bundle, where the typical
fiber is a positive dimensional simple unital algebra, completely determines
the bundle structure. Furthermore, we derive a homomorphism version
of the Shanks-Pursell theorem and identify a class
of homomorphisms
of multiplicative semigroups between spaces of smooth functions
on finite dimensional manifolds,
including  all isomorphisms.
Our approach is based on a method called range decreasing
group homomorphisms.
\end{abstract}


\maketitle

\pagestyle{myheadings} \markboth{\centerline{N. Zhang}}{\centerline{Range decreasing group homomorphisms}}


\section{Introduction \label{intro}}

In this paper, we assume that all manifolds and topological spaces are Hausdorff, and that all locally convex spaces are sequentially complete. Unless otherwise specified, we denote by $V$ (respectively $W$) a finite dimensional $C^{\infty}$ manifold,
possibly with boundary, that may not be paracompact or connected, and by $G$ a finite or infinite dimensional locally exponential connected Lie group modelled on a locally convex space (notably, all Banach Lie groups are locally exponential).
Let $\cf$ be a fixed smoothness class among $C^k_{MB}$ ($C^k$ in the Michal-Bastiani sense),
$k=0, 1, \cdots, \infty$, and let $\pi_{\scg, V}: \scg \to V$ be a $C^{\infty}$ Lie group
bundle with typical fiber $G$. The space $\fvscg$ of $\cf$ sections of the bundle $\pi_{\scg, V}$
is a group under pointwise group operation.
If $G$ is a Banach Lie group, we may also take $\cf$ to be
$C^k$, $k=0, 1, \cdots, \infty$, in the sense of Fr\'echet
differentiability.
In the case where $G$ is finite dimensional, we allow $\cf$ to be
locally H\"older $\ka$, $k=0, 1, \cdots$, $0 \le \alpha \le 1$,
where $C^{k, 0}=C^k$, or locally Sobolev $\kp$, $k=1, 2, \cdots$, $1 \le p <\infty$, $kp>\dimv$, or $k = \dimv$ and $p=1$.
The groups $\fvscg$ are useful. For example,
the gauge group of a principal $G$-bundle is isomorphic to the group of $C^{\infty}$ sections of
the associated Lie group bundle with typical fiber $G$, defined by the
conjugation action of $G$ on itself.

Suppose that $\pi_{\scg_1, V}: \scg_1 \to V$ and $\pi_{\scg_2, W}: \scg_2 \to W$ are $C^{\infty}$ Lie group
bundles with the same typical fiber $G$ and $\cf$, $\tilde{\cf}$ are
two smoothness classes. Let $H$ be a subgroup of $\fvscgo$, and let
$f: H \to \ftscgt$ be a group homomorphism.
We define $f$ as a weighted composition operator
if there exist a map $\phi: W \to V$ and a (not necessarily continuous) bundle morphism
   $\gamma: \phi^{\ast} \scg_1 \to \scg_2$ over the identity map on $W$
   such that the restriction of $\gamma$ to each fiber is a group isomorphism
   $(\phi^{\ast} \scg_1)_w \to (\scg_2)_w$ and
    \begin{equation} \label{weightedco}
      f(x)=\gamma \circ (\phi^{\ast} x), \hspace{2mm} x \in H,
    \end{equation}
    where $\phi^{\ast} \scg_1$ (respectively $\phi^{\ast} x$) is the total space
    of the pullback of the bundle $\pi_{\scg_1, V}$
    (respectively the pullback of the section $x$) by $\phi$ (even if $\phi$ is not continuous,
    $\phi^{\ast} \scg_1$ as a set (respectively $\phi^{\ast} x$
    as a map $W \to \phi^{\ast} \scg_1$) is still well-defined).
    If the Lie group bundles $\pi_{\scg_1, V}$ and $\pi_{\scg_2, W}$ are vector bundles (respectively algebra bundles)
    and $f$ is a linear map (respectively an algebra homomorphism),
    we further require that the restriction of $\gamma$ to each fiber
    is a linear isomorphism (respectively an algebra isomorphism).
 When the Lie group bundles $\pi_{\scg_1, V}$ and $\pi_{\scg_2, W}$ are trivial,
 the space $\fvscgo$ (respectively $\ftscgt$)
 corresponds to the mapping group $\fvg$ (respectively $\tilde{\cf}(W, G)$).
 In this scenario, the bundle morphism $\gamma$ can be interpreted as a map $W \times G \to G$,
 where $\gamma(w, \cdot): G \to G$ is a group isomorphism for each $w \in W$.
 Consequently, the relation (\ref{weightedco}) takes the form
 \begin{equation} \label{wcotrivial}
 f(x)(w)=\gamma(w, x \circ \phi(w)), \hspace{2mm} w \in W, \hspace{2mm} x \in H.
 \end{equation}
 If each $\gamma(w, \cdot)$ is the identity map on $G$, we refer to $f$ as a composition or pullback operator.

 Weighted composition operators have been the focus of research for many years, with various special cases explored extensively.
   Let $\bff=\br$ or $\bc$,
$X$ a completely regular topological space and $\ca$ a subalgebra of $C(X, \bff)$. We say that
$X$ is $\ca$-realcompact if every nonzero algebra homomorphism
$f: \ca \to \bff$ is an evaluation map $E_a: \ca \ni x \mapsto x(a) \in \bff$ for some point $a \in X$.
When $\ca=C(X, \br)$,  this condition reduces to the
usual realcompactness, originally introduced by Hewitt in \cite{he} under the term “Q-space”.
Suppose $X$ is $\ca$-realcompact, $S$ is a nonempty set and $\cm(S, \bff)$
is the space of all maps $S \to \bff$. It is clear that any algebra homomorphism $f: \ca \to \cm(S, \bff)$
such that $E_s \circ f \not\equiv 0$ for each $s \in S$
is the pullback by a map $\phi: S \to X$.
We say that a manifold $V$ is smoothly realcompact, if it is $C^{\infty}(V, \bff)$-realcompact \cite[Chapter IV]{kmbook}.
In \cite{g, mr}, Grabowski and Mr\v{c}un demonstrate that every algebra isomorphism $C^{\infty}(V, \bff) \to C^{\infty}(W, \bff)$,
where $V$ and $W$ are positive dimensional manifolds which are not necessarily paracompact
or smoothly realcompact, is the pullback by a $C^{\infty}$ diffeomorphism $W \to V$. This result generalizes earlier findings \cite{gk, myers, pursell}.
Any algebra isomorphism $C^k(V, \bff) \to C^k(W, \bff)$ is a linear biseparating map.
A theorem of Araujo states that any linear biseparating map between spaces of
$C^k$ (where $k \in \bn$) functions, defined on open subsets of positive dimensional Euclidean spaces
and taking values in Banach spaces, must be a weighted composition operator
\cite{ar}.

Let $\xv$ be the Lie algebra of all $C^{\infty}$ vector fields with compact support on $V$.
The Shanks-Pursell Theorem asserts that if there exists a Lie algebra isomorphism $f:$ $\xv$ $\to$ $\xw$,
then the manifolds $V$ and $W$ are diffeomorphic \cite{sp}.
Furthermore, there exists a $C^{\infty}$ diffeomorphism $\psi: V \to W$ such that $f = d\psi$ \cite{kmo}.
Consequently, $f$ is a weighted composition operator.
For further Shanks-Pursell type results concerning isomorphisms of Lie algebras of vector fields or
isomorphisms of Lie algebras of linear differential operators, see \cite[Section X]{o} and \cite{gim, ab, gp}.

According to a result by Milgram \cite[Theorem A]{mi}, any isomorphism $C(X, \br) \to C(X', \br)$
of multiplicative semigroups, where $X$ and $X'$ are compact topological spaces,
can be expressed as a weighted composition operator  as described
in (\ref{wcotrivial}).
In this case, $\phi: X' \to X$ is a homeomorphism and
each map $\gamma(a', \cdot): \br \to \br$ for $a' \in X'$
is an isomorphism of multiplicative semigroups.
Furthermore, Mr\v{c}un and \v{S}emrl establish that
any isomorphism $C^k(V, \mathbb{R}) \to C^k(W, \mathbb{R})$,
$k \in \bn$, of multiplicative semigroups
is the pullback by a $C^k$ diffeomorphism \cite{ms}.

Sufficient conditions for a homomorphism $f$ between groups of continuous maps
taking values in a topological group $G$
to be a weighted composition
operator are provided
in \cite{hr, fe}. In particular,
homomorphisms $f$ which commute with  continuous endomorphisms of
the group $G$ are considered by Hern\'{a}ndez and R\'{o}denas in \cite{hr}; and weakly biseparating
isomorphisms $f$ are investigated by Ferrer, Gary and Hern\'{a}ndez in \cite{fe}.
Additionally, a special case of the weighted composition operator
as defined in (\ref{wcotrivial}) is explored in \cite{z23},
where $V$ is compact and $G$ is finite dimensional.

It turns out that the appropriate domain for the weighted composition operator $ f $
described in (\ref{weightedco}) is the following subgroup of $\fvscgo$.
Let $\one=\one_{\pi_{\scg_1, V}} \in \fvscgo$ be the section that
assigns the identity element of the fiber $(\scg_1)_v$ to each point $v \in V$.
We say that the support of a section $x \in \fvscgo$ is compact
if $x$ coincides with $\one$ outside a compact subset of $V$.
We denote by
 $\fcvscgo \subset \fvscgo$ the subgroup of sections with compact support.
 Additionally, we define
 $\fczvscgo \subset \fcvscgo$ \label{fczvscgo} as the subgroup of sections $x$ for which
there exist a compact subset $K \subset V$ with $\supp x \subset K$
  and a homotopy $H: [0, 1] \times V \to \scg_1$
  relative to $V \setminus K$ such that $H(0, \cdot)=\one$, $H(1, \cdot)=x$ and
  $H(t, \cdot) \in \fvscgo$ for all $t \in [0, 1]$.
  If the bundle $\pi_{\scg_1, V}$ is trivial,
  we write $\fczvg$ for $\fczvscgo$.
  When $V$ is compact, the space $\fvscgo$ is a Lie group, and $\fczvscgo$
  represents the component of $\fvscgo$ containing the identity element $\one$.
  In the case where $\pi_{\scg_1, V}$ is a vector bundle $\pi_{\xi, V}$, we have that $\cf_c^0(\pi_{\xi, V}) = \cf_c(\pi_{\xi, V})$.
  There exist group homomorphisms $f_0: \fvscgo \to \ftscgt$ whose restrictions to $\fczvscgo$
  are weighted composition operators.
  However, the homomorphisms $ f_0 $ are not weighted composition operators
  on $\fvscgo$.

  This paper presents a necessary and sufficient condition for
  a group homomorphism $f: \fczvscgo \to \ftscgt$ to be a weighted composition
  operator (Theorem \ref{lgb}).
  Furthermore, the above result can be significantly simplified under specific conditions: if the bundles $\pi_{\scg_1, V}$ and $\pi_{\scg_2, W}$ are trivial (discussed in Section \ref{mappinggroup}), or if both bundles are vector bundles of the same rank $r \in \mathbb{N}$ and $f$ is linear (Theorem \ref{finiterankbundle}).

Here are several applications of our results.

Application 1: Suppose that $A_1$ (respectively $A_2$) is a positive dimensional simple
unital $\bff$-algebra (all algebras in this paper are associative), where $\bff=\br$ or $\bc$,
and $\pi_{\sa_1, V}: \sa_1 \to V$ (respectively $\pi_{\sa_2, W}: \sa_2 \to W$) is a $C^{\infty}$ $\bff$-algebra bundle
with the typical fiber $A_1$ (respectively $A_2$).
There is a wealth of nontrivial bundles of this type.
We prove that if there exists an algebra isomorphism
$C^{\infty}_c(\pi_{\sa_1, V})$ $\to$ $C^{\infty}_c(\pi_{\sa_2, W})$ (respectively $C^{\infty}(\pi_{\sa_1, V})$ $\to$
$C^{\infty}(\pi_{\sa_2, W})$),
  then the algebra bundles $\pi_{\sa_1, V}$ and $\pi_{\sa_2, W}$ are isomorphic (Corollary \ref{isoalgebra}).
  Hence
the algebraic structure of $C^{\infty}_c(\pi_{\sa_1, V})$ (respectively $C^{\infty}(\pi_{\sa_1, V})$) completely
determines the bundle structure of $\pi_{\sa_1, V}$.
Moreover, if $A_1=A_2$, then any $\bff$-algebra homomorphism $f: \fcsao \to \ftsat$
(respectively $f: \fsao \to \ftsat$) for which
$E_{w} \circ f|_{\fcsao} \not\equiv \zero$, $w \in W$, must be a weighted composition operator (Theorem \ref{finitesp}).
Thus every positive dimensional manifold $V$ is $\fcvf$-realcompact,
even though it is not necessarily realcompact.

Application 2: We obtain a homomorphism version of the Shanks-Pursell Theorem (Theorem \ref{shanks}),
in contrast to the traditional isomorphism.

Application 3: We provide
sufficient conditions for a homomorphism $\fvr$ $\to$ $\tfwr$
of multiplicative semigroups to be a weighted composition operator as in (\ref{wcotrivial}),
where each $\gamma(w, \cdot): \br \to \br$ for $w \in W$
is an isomorphism of multiplicative semigroups
(Theorem \ref{smgp}).
As a consequence, we extend Milgram's result to the context
 of semigroup isomorphisms $C(V, \br) \to C(W, \br)$, where
 $V$ and $W$ are arbitrary finite dimensional manifolds.
Additionally, we generalize the theorem of Mr\v{c}un and \v{S}emrl to accommodate
 the case when $k=\infty$
(Corollary \ref{semigpiso}).

Application 4: We demonstrate that any linear bijection $f:$ $\fvrn$ $\to$ $\fwrn$, $n \in \bn$,
that satisfies
  $f( \cf(V, \br^n \setminus \{\zero \}))$ $=$ $\cf(W, \br^n \setminus \{\zero\})$
must be a weighted composition operator (Theorem \ref{diffeo}).

Our approach utilizes the technique of range decreasing group homomorphisms,
originally introduced in \cite{z23} and further explored in this paper.
Let $S_j$, $j=1, 2, 3$, be nonempty sets, $\cm(S_i, S_3)$ the space of all maps
$S_i \to S_3$, $i=1, 2$, and $\mathscr{X}$ a subset of $\cm(S_1, S_3)$.
We say that a map $f: \mathscr{X} \to \cm(S_2, S_3)$ is range decreasing, if
$$f(x)(S_2) \subset x(S_1), \hspace{2mm} x \in \mathscr{X}.$$
Any composition operator defined by $f(x)=x \circ \phi$,
where  $x \in \cm(S_1, S_3)$ and $\phi$ is a map $S_2 \to S_1$, is inherently range decreasing.
On the other hand,
a linear range decreasing map $\fvr \to \tilde{\cf}(W, \br)$ is not always
a composition operator \cite[P. 2190]{z23}. Furthermore,
for any $n \in \bn$,
there exist $C^{\infty}$ range decreasing maps from $C^{\infty}(S^2, \bc^n)$ to itself
that are not composition operators \cite[P. 2180]{z23}.
The work presented in \cite{z23} uncovers unexpected findings. For instance,
if $V$ and $W$ are compact manifolds, possibly with boundary,
and $G$ is a finite dimensional connected Lie group with $\dim_{\br} G \ge 2$
(in this case, $\fvg$ is a Banach or Fr\'echet Lie group), then any range decreasing group
homomorphism $f: \fvg \to \tilde{\cf}(W, G)$ must be a composition
operator on any component of $\fvg$ containing an element $x_0$ such that $x_0(V)$
is nowhere dense in $G$ \cite[Theorem 1.1]{z23}.

In this paper, we broaden the scope of previously established results to encompass much more general contexts.
Let $V$, $W$ and $G$ be as in the first paragraph
of this section, and let $f: \fcvg \to \sg$ (respectively $f: \fvg \to \sg$)
be a range decreasing group homomorphism
  with $E_s \circ f|_{\fcvg} \not\equiv \one$ for each $s \in S$.
  If $\dim_{\br} G \ge 2$,
  then we can find a map $\phi: S \to V$ such that $f(x)=x \circ \phi$ for
  any $x \in \fczvg$ and for any $x$ in the domain of $f$ with $x(V) \not=G$.
  Furthermore, if $ G $ is a locally convex space, this relationship holds for all
  $ x $ in the domain of $ f $ (Sections \ref{additive} and \ref{rdcomposition}).
  We also construct a range decreasing group homomorphism from $C^{\infty}(S^3, \mathrm{SU}(2))$ to itself
that is not a composition operator for certain components  of $C^{\infty}(S^3, \mathrm{SU}(2))$
  (Proposition \ref{su2}).
  Since $C^{\infty}$
partition of unity is not always available on $\supp x$, where $x \in \fvg$,
the proofs for the case of $f: \fvg \to \sg$ differ significantly from those presented in \cite{z23}.

The method of range decreasing group homomorphism is independent of the topology of the spaces under consideration.
This approach is applicable to group homomorphisms in general, rather than being limited to group isomorphisms or linear maps.
It can be applied to any finite dimensional manifold, regardless of its paracompactness or realcompactness,
and is suitable for a wide range of smoothness classes.
Additionally, the “range decreasing” condition can be combined with holomorphicity,
rather than group homomorphism, in the exploration of composition operators.
For instance,
it is proved in \cite{z23} that oftentimes a range decreasing holomorphic
map $\cf(V, M)$ $\to$ $\tilde{\cf}(W, M)$, where $V$, $W$ are compact manifolds
and $M$ is a finite dimensional complex manifold without boundary,
is a composition operator (for the complex structure on $\fvm$, see \cite{l04, ls}).

This paper is organized as follows.  In Section 2, we review relevant
facts regarding mapping spaces and prove several propositions which will be
necessary for subsequent discussions.
In Sections 3, 4 and 5, we study range decreasing group homomorphisms.
In Section \ref{grouph}, we prove that if $\dim_{\br} G \ge 2$,
any nonconstant range decreasing group homomorphism $\fczvg \to G$
is the evaluation $E_{\vb}$ at some point $\vb$ of $V$.
Let $\sce$ be a locally convex space with $\dim_{\br} \sce \ge 2$. In Section \ref{additive},
we show that any range decreasing (additive) group homomorphism $\fve \to \sce$ with $f|_{\fcve}$ $\not\equiv$ $\zero$
takes the form $E_{\vb}$ on $\fve$.
In Section \ref{rdcomposition}, we establish that if $\dim_{\br} G \ge 2$,
any range decreasing group homomorphism
$f: \fcvg \to G$ (respectively $f: \fvg \to G$) with $f|_{\fczvg} \not\equiv \one$
can be expressed as $E_{\vb}$ on the subset of $\fcvg$ (respectively of $\fvg$)
consisting of non-surjective maps.

Let $f: \fczvscgo \to \ftscgt$ be a group homomorphism that can be represented
as a weighted composition operator, as described in (\ref{weightedco}).
Our main results in Section \ref{wco1} demonstrate that, under certain mild conditions,
the bundle morphism $\gamma$ is automatically continuous and  the map $\phi$ is automatically smooth.
However, in the absence of any conditions,
 $\gamma$ may be discontinuous.

Sections \ref{mappinggroup} to \ref{ah} highlight the applications of the results discussed in
the earlier sections.
In these sections, necessary and sufficient conditions,  along with sufficient conditions,
are provided for various types of homomorphisms,
including group homomorphisms, linear maps, semigroup homomorphisms, Lie algebra homomorphisms and algebra homomorphisms,
to be classified as weighted composition operators.
In Section \ref{mappinggroup},
we focus on group homomorphisms of the form
$f: \fvg \to \tilde{\cf}(W, G)$.
In Section \ref{lm}, we explore Application 4, which concerns linear maps.
Section \ref{hsg} presents an analysis
 of semigroup homomorphisms, leading to the results in
 Application 3.
Section \ref{hlgb}
addresses group homomorphisms $f: \fczvscgo \to \ftscgt$.
In Section \ref{vf}, we offer a generalization of the Shanks-Pursell Theorem, as outlined in Application 2.
Finally, in Section \ref{ah}, we classify specific types of algebra homomorphisms and conclude with Application 1.


\section{Preliminaries \label{pre}}

The definition of a $C^k$ map, where $k=1, 2, \cdots, \infty$, in the sense of Fr\'echet
differentiability
from an open subset of a Banach space to another Banach space
 is based on the concept of
differentials \cite[Section 2.3]{amr},
while a $C^k_{MB}$ map ($C^k$ in the Michal-Bastiani sense)
from an open subset of a locally convex space to another
locally convex space is defined in terms of directional
derivatives \cite[Section I.2]{towards}.
For a map between open subsets of Banach spaces, $C^{k+1}_{MB}$ implies $C^k$ in the sense of Fr\'echet differentiability, which in turn implies $C^k_{MB}$. Therefore $C^{\infty}_{MB}=C^{\infty}$.
Both types of smoothness classes can be generalized to maps between manifolds.

A locally convex Lie group $G$ is called locally exponential
if it has a $C^{\infty}$ exponential map $\exp _G: \fg \to G$, where $\fg$ is the Lie algebra of $G$,
such that there exists an open neighborhood $U$ of $\zero \in \fg$ that is
mapped diffeomorphically onto an open neighborhood of $\one \in G$
\cite[Section IV.1]{towards}.
Let $S$ be a nonempty set. The constants form a subgroup of the mapping group $\msg$
which we identify with $G$.

\begin{prop} \label{algebrard}
  Let $f: \fvg \to \msg$ be a group homomorphism such that $f(a)=a$ for each
  $a \in G$ and
  $f\left( \cf(V, G \setminus \{\one\}) \right) \subset \cm(S, G \setminus \{\one\})$.
  Then $f$ is range decreasing.
\end{prop}
\begin{proof}
  We may assume that $S$ consists of a single point (i.e. $\cm(S, G)=G$). Otherwise
  consider $E_s \circ f$ for each $s \in S$.
  Note that $f(x f(x)^{-1})=\one$, $x \in \fvg$.
  Thus $\one \in (x f(x)^{-1})(V)$, which implies that $f(x) \in x(V)$.
\end{proof}

Let $G$ (respectively $A$) be a locally exponential connected locally
convex Lie group (respectively a positive dimensional unital $\bff$-algebra, where $\bff=\br$ or $\bc$).
A $C^{\infty}$ Lie
group bundle $\pi_{\scg, V}: \scg \to V$ (respectively an $\bff$-algebra bundle \label{abundle} $\pi_{\sa, V}: \sa \to V$)
with typical fiber $G$  (respectively $A$) is defined as a $C^{\infty}$ fiber bundle (respectively vector bundle)
such that for any $v \in V$, the fiber $\scg_v$ (respectively $\sa_v$) is equipped with a Lie group (respectively an algebra) structure and
there are an open neighborhood $O$ of $v$ and a diffeomorphism $O \times G \rightarrow
\pi^{-1}_{\scg, V}(O)$ (respectively $O \times A \rightarrow
\pi^{-1}_{\sa, V}(O)$) preserving the Lie group (respectively the algebra) structure fiberwisely.
Any Lie group bundle $\pi_{\scg, V}: \scg \to V$ induces a Lie algebra bundle $\hat{\pi}_{\xi, V}: \xi \to V$
whose fiber $\xi_v$ at each $v \in V$ corresponds to the Lie algebra of $\scg_v$.
We write $\exp_{\scg}: \xi \to \scg$ for the $C^{\infty}$ exponential map, which can be defined fiberwisely.
For more information about Lie group bundles, see \cite{lo}.

  \begin{prop} \label{relativeh}
  Let $\Omega \subset \xi$ be an open neighborhood of the zero section of the bundle $\hat{\pi}_{\xi, V}$
such that the restriction of the exponential map $\exp_{\scg}|_{\Omega}: \Omega \to \exp_{\scg}(\Omega)$ is a $C^{\infty}$ diffeomorphism.
If there exist a precompact open subset $O \subset V$
and a homotopy $H: [0, 1] \times V \to \scg$
relative to $V \setminus O$ such that  $x_t=H(t, \cdot) \in \fvscg$ for each $t \in [0, 1]$,
then there exist finitely many $y_1, \cdots, y_p \in \fcvscg$ such that $\supp y_j \subset \overline{O}$,
$y_j(V) \subset \exp_{\scg}(\Omega)$, $j=1, \cdots, p$, and
$x_0 x_1^{-1}=y_1 \cdots y_p.$
\end{prop}
\begin{proof}
  Due to the compactness of $[0, 1] \times \overline{O}$, there exists a finite sequence
of real numbers
$0=t_0<t_1<\cdots<t_p=1$
such that
$x_{t_{j-1}} x_{t_{j}}^{-1} \in \cf_c(V, \exp_{\scg}(\Omega))$ for $j=1, \cdots, p$.
Set $y_j=x_{t_{j-1}} x_{t_{j}}^{-1}$.
\end{proof}

 Recall the subgroup $\fczvscg$ of $\fcvscg$ as defined in Section \ref{intro}.
  For any $\xh \in \cf_c(\hat{\pi}_{\xi, V})$, we have $\exp_{\scg} \circ \xh \in \fczvscg$. By Proposition \ref{relativeh},
 the subgroup $\fczvscg$ \label{generate} is generated by the subset $ \{\exp_{\scg} \circ \xh: \xh \in \fcvx \}$.

 If $V$ is compact, then the mapping group $\fvg$
carries a natural Lie group structure.
Let $\exp_G: \fg \to G$ be the exponential map of $G$.
The Lie algebra associated with
 $\fvg$ is identified with $\cf(V, \fg)$, and
the exponential map of $\fvg$ is given by
$\cf(V, \fg) \ni \xh \mapsto \exp_G  \circ \xh \in \fvg$, see \cite[Theorem II.2.8]{towards} and \cite[Section 4(G)]{kr}.
Since $G$ is locally exponential, it follows
that $\fvg$ is also locally exponential.
Moreover, $\fczvg$ represents the component of the Lie group $\fvg$
that contains the identity element $\one$.

\begin{prop} \label{phi}
  Let $M$ be a $C^{\infty}$ connected locally convex manifold without boundary, where $\dim_{\br} M>0$.
 If there exists a map $\phi: W \to V$ such that $x \circ \phi \in \tilde{\cf}(W, M)$ for any
$x \in \fvm$, then $\phi$ is an $\tilde{\cf}$ map. Furthermore, if $M$
is a Lie group $G$, the same conclusion holds under a less stringent condition. Specifically, it suffices to assume that
$x \circ \phi \in \tilde{\cf}(W, G)$ for any $x \in \fczvg$ with $x(V) \not=G$.
\end{prop}
\begin{proof}
  Suppose that $M$ is modelled on the locally convex space $\sce$,
  and consider a local chart $(U, \Phi)$ of $M$ such that $\Phi(U)$
  is an open neighborhood of $\zero \in \sce$. If we assume that $M$ is a Lie group $G$,
  then $\sce$ represents the Lie algebra of $G$.
  In this context, we further specify that $U$ is a small neighborhood of $\one \in G$,
  and that the chart is given by $\Phi=\exp^{-1}_G|_U$, where $\exp_G$ is the exponential map of $G$.
  Take
  $$a \in \Phi(U) \setminus \{\zero\} \hspace{1.5mm} \text{with} \hspace{1.5mm}
  ta \in \Phi(U), \hspace{1.5mm}  -2 \le t \le 2.$$
  Let $w_0 \in W$, $(\fv_0, \Psi)$ a local chart of $V$
  with $\phi(w_0) \in \fv_0$, $\ve>0$, $\chi_0 \in C^{\infty}_c(V, \br)$
   and $\fv_1$ a precompact open neighborhood
  of $\phi(w_0)$ such that $0 \le \chi_0(v)<\chi_0(\phi(w_0))$ for any $v \in V \setminus \{\phi(w_0)\}$
  and $\overline{\fv_1} \subset \fv_0$. Furthermore, we require that the condition
  $|\chi_0(v)-\chi_0(\phi(w_0))|/\chi_0(\phi(w_0))<\ve$ holds if and only
  if $v \in \fv_1$. Set
  $$x_0=\Phi^{-1} \circ \left(\chi_0 a/\chi_0(\phi(w_0)) \right) \in C^{\infty}(V, M).$$
  If $M=G$, then $x_0 \in \fczvg$ and $x_0(V) \not=G$.
  Since $x_0 \circ \phi$ is continuous, we can find an open neighborhood $\fw_0$ of $w_0$
  with
  $\phi(\fw_0) \subset \fv_1.$

  Let $\Psi_1, \cdots, \Psi_{\dim_{\br} V}: \fv_0 \to \br$ be the components
  of $\Psi$. Take $\chi_1 \in C^{\infty}_c(\fv_0, \br)$ with $\chi_1|_{\fv_1} \equiv 1$.
  Define
  $$x_j=\Phi^{-1} \circ \left( \chi_1 \Psi_j a/||\chi_1 \Psi_j||_C  \right) \in C^{\infty}(V, M),
  \hspace{2mm} j=1, \cdots, \dim_{\br} V,$$
  where $||\chi_1 \Psi_j||_C=\sup_{v \in \fv_0}|\chi_1 \Psi_j(v)|$.
  If $M=G$, then  $x_j \in \fczvg$ and $x_j(V) \not=G$.
  Note that $x_j \circ \phi \in \tilde{\cf}(W, M)$.
  By the equation above, we have
  $$(\Psi_j \circ \phi)(w) a=||\chi_1 \Psi_j||_C (\Phi \circ x_j \circ \phi)(w), \hspace{1.5mm} w \in \fw_0, \hspace{1.5mm} j=1, \cdots, \dim_{\br} V.$$
  Hence $\phi$ is an $\tilde{\cf}$ map.
\end{proof}

\begin{prop} \label{measurable}
  Let $(X, \mathfrak{M})$ be a measurable space, $(Y, d)$ a metric space, $\cb_Y$
the Borel $\sigma$-algebra of $Y$ and $\psi_0: X \to Y$ a map. If there exists
a sequence of $(\mathfrak{M}, \cb_Y)$-measurable maps $\{\psi_j: X \to Y \}$
 such that
$\psi_0(\alpha)=\lim_{j \to \infty} \psi_j(\alpha)$ for every $\alpha \in X$,
then $\psi_0$ is $(\mathfrak{M}, \cb_Y)$-measurable.
\end{prop}
\begin{proof}
  It suffices to show that $\psi_0^{-1}(O) \in \mathfrak{M}$ for any open subset $O \subset Y$.
  If $O \not=Y$, we define
$O_m=\{\beta \in O : d(\beta, Y \setminus O)>1 /(m+1)\}$, where $m \in \mathbb{N}$.
The set $O_m$ is open in $Y$. Note that $\alpha \in \psi_0^{-1}(O)$ if and only if there exist
positive integers $m$ and $N$
such that $\psi_j(\alpha) \in O_m$ for each $j \ge N$.
So
$\psi_0^{-1}(O)=\bigcup_{m \in \mathbb{N}} \bigcup_{N \in \mathbb{N}} \bigcap_{j \geq N} \psi_j^{-1}\left(O_m\right) \in \mathfrak{M}.
$
\end{proof}


\section{Range decreasing homomorphisms $\fcvg \to G$ \label{grouph}}

The main result of this section is Theorem \ref{nbhd1}.

Let $\pi_{\scg, V}: \scg \to V$ be a $C^{\infty}$ Lie group bundle with typical
fiber $G$, $\hat{\pi}_{\xi, V}: \xi \to V$ the Lie algebra bundle of $\pi_{\scg, V}$,
$\exp_{\scg}: \xi \to \scg$ the exponential map of $\pi_{\scg, V}$,
$H$ a group and $\varphi: \fczvscg \to H$ a nonconstant group homomorphism.
We say that $v_0 \in V$ is a support point
of $\varphi$ if, for every open neighborhood $O$ of $v_0$, there exists a section $\xh \in \cf_c(\hat{\pi}_{\xi, V})$
such that $\supp \xh \subset O$ and $\varphi(\exp_{\scg} \circ \xh) \not=\one \in H$. Put
$$\cs_{\varphi}=\{v_0 \in V: v_0 \hspace{1mm} \text{is a support point of} \hspace{1mm} \varphi\}.$$

\begin{prop} \label{sp}
  The subset $\cs_{\varphi}$ of $V$ is nonempty and closed.
  Furthermore, if $\xh_1, \xh_2 \in \cf_c(\hat{\pi}_{\xi, V})$ have the same germ
at each point of $\cs_{\varphi}$, then
$\varphi(\exp_{\scg} \circ \xh_1)=\varphi(\exp_{\scg} \circ \xh_2)$.
\end{prop}
\begin{proof}
It follows from the definition of $\cs_{\varphi}$ that $V \setminus \cs_{\varphi}$ is open.
Suppose that $\xh_1, \xh_2 \in \cf_c(\hat{\pi}_{\xi, V})$ have the same germ
at each point of $\cs_{\varphi}$. Let $v \in \supp \xh_1 \cup \supp \xh_2$.
If $v \in \cs_{\varphi}$, we can select a precompact open neighborhood $O_v$ of $v$ such that
$\xh_1|_{O_v}=\xh_2|_{O_v}$. Conversely, if $v \not\in \cs_{\varphi}$, we can choose a precompact
open neighborhood $O_v \subset V \setminus \cs_{\varphi}$ of $v$ such that $\varphi(\exp_{\scg} \circ \xh)=\one$
for every $\xh \in \cf_c(\hat{\pi}_{\xi, V})$ with $\supp \xh \subset O_v$.
Take a finite subcover
  $\{O_{v_j}\}_{j=1}^q$ of the open cover $\{O_v\}$ of $\supp \xh_1 \cup \supp \xh_2$
  and a $C^{\infty}$ partition of unity $\{\chi_j \}_{j=1}^q$ subordinate to $\{O_{v_j}\}_{j=1}^q$.
  Note that
  $[\chi_{j_1} \xh_i, \chi_{j_2} \xh_i] \equiv \zero$, $j_1, j_2=1, \cdots, q$, $i=1, 2$,
  where $[\cdot, \cdot]$ is the Lie bracket on $\cf_c(\hat{\pi}_{\xi, V})$,
  which can be defined fiberwisely. Hence
  \begin{equation*}
  \te \exp_{\scg} \circ \xh_i=\exp_{\scg} \circ (\Sigma_{j=1}^q \chi_j \xh_i)=\Pi_{j=1}^q \exp_{\scg} \circ (\chi_j \xh_i), \hspace{2mm} i=1, 2.
  \end{equation*}
  This implies that
  $\varphi(\exp_{\scg} \circ \xh_1)=\varphi(\exp_{\scg} \circ \xh_2)$.
  If $\cs_{\varphi}=\emptyset$, then $\varphi(\exp_{\scg} \circ \xh_i)$ $=$ $\one$ for $i=1, 2$.
  Thus $\varphi$ is constant, which leads to a contradiction.
\end{proof}

\begin{prop} \label{phicontinuous}
  Let $\pi_{\scg_1, V}: \scg_1 \to V$ and $\pi_{\scg_2, W}: \scg_2 \to W$ be $C^{\infty}$ Lie group
bundles (where the typical fiber of $\pi_{\scg_1, V}$ is not necessarily the same as that of $\pi_{\scg_2, W}$),
 and let $f: \fczvscgo \to C(\pi_{\scg_2, W})$ be a group homomorphism
  such that for any $w \in W$, the set $\cs_{E_w \circ f}$ consists of
  a single point $v(w)$. Then the map $\phi: W \ni w \mapsto v(w) \in V$
  is continuous.
\end{prop}
\begin{proof}
    Let $\hat{\pi}_{\xi_1, V}: \xi_1 \to V$ be the Lie algebra bundle of $\pi_{\scg_1, V}$,
     and let $\exp_{\scg_1}: \xi_1 \to \scg_1$ be the exponential map of $\pi_{\scg_1, V}$.
    For any $w_0 \in W$ and any open neighborhood $O$ of $\phi(w_0)$,
    take $\xh_0 \in \cf_c(\hat{\pi}_{\xi_1, V})$
    such that $\supp \xh_0 \subset O$ and $f(\exp_{\scg_1} \circ \xh_0)(w_0)\not=\one  \in (\scg_2)_{w_0}$.
    If $W_0$ is a neighborhood
    of $w_0$ such that $f(\exp_{\scg_1} \circ \xh_0)(w) \not=\one  \in (\scg_2)_{w}$, $w \in W_0$,
    then by Proposition \ref{sp}, $\phi(W_0) \subset \supp \xh_0 \subset O$.
    Hence $\phi$ is continuous at $w_0$.
\end{proof}

\begin{lemma} \label{germ}
  Let $\sce$ be a locally convex space, and let $f: \fcve \to \sce$ be
  a nonzero range decreasing (additive) group homomorphism.
  If $\cs_{f}$ consists of a single point $\vb$, then $f$
  is the evaluation $E_{\vb}$ at $\vb$.
\end{lemma}
\begin{proof}
It follows from Proposition \ref{sp} that for any $x \in \fcve$,
the value $f(x)$ only depends on the germ of $x$ at $\vb$.
Assume for contradiction that there exists $x_1 \in \fcve$ such that $f(x_1) \not= x_1(\vb)$.
Choose $C^{\infty}$ functions $\mu: \sce \to \br$
and $\chi: \br \to [0, 1]$ such that
$\mu\left(f(x_1)\right)=1$, $\mu\left(x_1(\vb)\right)=0,$ $\supp \chi \subset \left(-1/2, 1/2 \right)$
and $\chi$ is constant $1$ on a neighborhood of $0 \in \br$. Set
$x_2=(\chi \circ \mu \circ x_1) x_1 \in \fcve$.
Then $x_1$ and $x_2$ have the same germ at $\vb$ and  $\mu \circ x_2(v)<1/2$ for every $v \in V$.
Hence $f(x_2)=f(x_1) \not\in x_2(V)$, this gives a contradiction.
\end{proof}

\begin{lemma} \label{rde}
  Let $f: \fcve \to \sce$ be a nonzero range decreasing group homomorphism. If $\dim_{\br} \sce \ge 2$,
  then there exists $\vb \in V$
  such that $f=E_{\vb}$.
\end{lemma}
\begin{proof}
We can express $\sce$ as a topological direct sum
$\sce_1 \oplus \sce_2$,
where $\sce_1$ and $\sce_2$ are closed subspaces of $\sce$ with $\dim_{\br} \sce_1, \dim_{\br} \sce_2 >0$.
Let $P_j: \sce \to \sce_j$ be the continuous projections for $j=1, 2$. It follows from
the range decreasing condition that $f(x_j) \in \sce_j$ for $x_j \in \cf_c(V, \sce_j)$.
Let $f_j=f|_{\cf_c(V, \sce_j)}$, $j=1, 2$.
Then
$f(x)=f_1 \circ P_1 \circ x+f_2 \circ P_2 \circ x$, $x \in \fcve$.

We claim that $f_j \not\equiv \zero$, $j=1, 2$. Note that $f$ is nonzero.
Let us assume that $f_1$ is nonzero.
Suppose $f_{2} \equiv \zero$.
Take $x_j \in \cf_c(V, \sce_j)$, $j=1, 2$,
such that $f_{1}(x_1) \not=\zero$ and $x_2$ does not vanish
on $\supp x_1$.
There exists $v_0 \in V$ such that
      $f_{1}(x_1) =f(x_1 +x_2 )=x_1(v_0) +x_2(v_0)$.
While the left hand side of the equation is in $\sce_1 \setminus \{\zero\}$, the right hand side is not.
We have a contradiction.
So $f_j \not\equiv \zero$ and $\cs_{f_j} \not= \emptyset$ (Proposition \ref{sp}), $j=1, 2$.

     Take $\vb \in \cs_{f_1}$. Next we show that
  $\cs_{f_1}=\cs_{f_2}=\{\vb\}.$
Assume for contradiction that there exists $v_2 \in \cs_{f_2}$ with $v_2 \not=\vb$. Take
an open neighborhood $\mathfrak{V}_1$ of $\vb$,
an open neighborhood $\mathfrak{V}_2$ of $v_2$ and
 $\xh_j \in \cf_c(V, \sce_j)$, $j=1, 2$, such that $\mathfrak{V}_1 \cap \mathfrak{V}_2=\emptyset$, $\supp \xh_1 \subset \mathfrak{V}_1$, $\supp \xh_2 \subset \mathfrak{V}_2$,
$f_1(\xh_1) \not=\zero$ and $f_2(\xh_2) \not=\zero$. Then
there exists $v_1 \in V$ such that
  $\xh_1(v_1) +\xh_2(v_1)=f_1(\xh_1)+f_2(\xh_2)$.
The left hand side of the equation is either in $\sce_1$ or
in $\sce_2$, but the right hand side is not contained in $\sce_1 \cup \sce_2$. We have
a contradiction. Therefore
$\cs_{f_2}=\{\vb\}$. Similarly $\cs_{f_1}=\{\vb\}$.
In view of Lemma \ref{germ}, we have $f=E_{\vb}$.
\end{proof}

Recall that each continuous homomorphism between locally exponential locally convex Lie groups
is $C^{\infty}$ \cite[Theorem IV.1.18]{towards}.

\begin{thm} \label{nbhd1}
  Let $f: \fczvg \to G$ be a nonconstant range decreasing group homomorphism
 and $\exp_G: \fg \to G$ the exponential map of $G$.
Then there exists a range decreasing linear map $\fh: \fcvfg \to \fg$ such that
\begin{equation} \label{fexp}
f(\exp_G \circ \xh)=\exp_G \circ \fh(\xh), \hspace{2mm} \xh \in \cf_c(V, \fg).
\end{equation}
  If $\dim_{\br} G \ge 2$, or $\dim_{\br} G=1$ and $\cs_{f}$ consists of a single point, then there exists $\vb \in V$
  such that
  \begin{equation} \label{fxvb}
  f(x)=x(\vb), \hspace{2mm}  x \in \fczvg.
  \end{equation}
\end{thm}
\begin{proof}
  Let $\cd$ be an open convex neighborhood of $\zero \in \fg$
such that $\exp_G|_{\cd}: \cd \to \exp_G(\cd)$ is a diffeomorphism.
Define a map
\begin{equation} \label{fh}
 \fh: \cf_c(V, \cd) \ni \xh \mapsto \left( \exp_G|_{\cd} \right)^{-1} \circ f(\exp_G \circ \xh) \in \cd.
\end{equation}
Since $f$ is range decreasing, $\fh$ is well-defined and $ \fh(\xh) \in \xh(V)$.

Let $\xh_1, \xh_2 \in \cf_c(V, \cd)$.
Take open neighborhoods $\tv_0$, $\tv$ of $\supp \xh_1 \cup \supp \xh_2$ such that
$\overline{\tv_0} \subset \tv$
and both $ \overline{\tv_0}$ and $\overline{\tv}$ are compact $C^{\infty}$ submanifolds of $V$, possibly with boundary.
Then $\cf(\overline{\tv}, G)$ is a Lie group modelled on the
locally convex space $\cf(\overline{\tv}, \fg)$. Set
\begin{eqnarray}
& \cf_{\overline{\tv_0}}(\overline{\tv}, G) = \{x \in \cf(\overline{\tv}, G): x|_{\overline{\tv} \setminus \tv_0} \equiv \one \} & \notag\\
&= \{x \in \cf(\overline{\tv}, G): x=y|_{\overline{\tv}}, \hspace{1mm}  y \in \fcvg, \hspace{1mm}  \supp y \subset \overline{\tv_0}\}. & \label{restriction}
\end{eqnarray}
If $\overline{\tv} \setminus \tv_0=\emptyset$, then $\tv$ is both
open and compact. Consequently, $\tv$
must be the union of finitely many compact components of $V$.
For any $x$ in the closed subgroup $\cf_{\overline{\tv_0}}(\overline{\tv}, G) \subset \cf(\overline{\tv}, G)$ (in infinite dimensions,
a closed subgroup is not always a Lie subgroup \cite[P. 374]{towards}),
take an open connected neighborhood $\su_x \subset \cf(\overline{\tv}, G)$
of $x$ such that $(x^{-1} y) (\overline{\tv}) \subset \exp_G(\cd)$ for each $y \in \su_x$. The maps
$$\su_x \ni y \mapsto (\exp_G|_{\cd})^{-1} \circ (x^{-1} y) \in \cf(\overline{\tv}, \cd) \subset \cf(\overline{\tv}, \fg)$$
provide local charts for $\cf(\overline{\tv}, G)$. It is clear that $\cf_{\overline{\tv_0}}(\overline{\tv}, G)$
is a Lie subgroup of $\cf(\overline{\tv}, G)$ modelled on the closed subspace of $\cf(\overline{\tv}, \fg)$
consisting of elements $\xh$ such that $ \xh|_{\overline{\tv} \setminus \tv_0} \equiv \zero$.
Let $ \cf^0_{\overline{\tv_0}}(\overline{\tv}, G)$
be the component of $ \cf_{\overline{\tv_0}}(\overline{\tv}, G)$  that contains  $\one$.
By (\ref{restriction}), $f$ induces a group homomorphism
$ f_1: \cf^0_{\overline{\tv_0}}(\overline{\tv}, G) \to G.$
It follows from the range decreasing condition that $f_1$ is continuous at $\one$.
So $f_1$ is $C^{\infty}$.
Let $\fh_1$ be the Lie algebra homomorphism of $f_1$. Note that
$\exp_G \circ \xh_1$, $\exp_G \circ \xh_2 \in \cf^0_{\overline{\tv_0}}(\overline{\tv}, G)$.
 According to (\ref{fh}),
\begin{equation} \label{f1}
\fh(\xh)=\fh_1(\xh), \hspace{2mm} \xh \in \cf_c(V, \cd) \cap \{t_1 \xh_1+t_2 \xh_2: (t_1, t_2) \in \br^2 \}.
\end{equation}
In particular, (\ref{f1}) is valid  for every $(t_1, t_2)$
in a sufficiently small neighborhood of $\zero \in \br^2$.

For any $\xh \in \fcvfg \setminus \cf_c(V, \cd)$, there exists $j \in \bn$ such that
$\xh/j \in \cf_c(V, \cd).$
Set $ \fh(\xh)=j\fh\left(\xh/j \right)$,
which is independent of $j$ ((\ref{f1})). Now
we obtain a well-defined linear map
 $\fh: \fcvfg \to \fg.$
It is clear that $\fh$ is range decreasing on $\fcvfg$.
By (\ref{fh}), we have (\ref{fexp}).
It follows from Lemmas \ref{germ} and \ref{rde}
that (\ref{fxvb}) holds for any $x = \exp_G \circ \xh$, $\xh \in \cf_c(V, \fg)$.
Hence it also holds for
any $x$ in the group $\fczvg$ generated by the subset $\{\exp_G \circ \xh: \xh \in \cf_c(V, \fg)\}$.
\end{proof}

If $V$ is  compact, then $\fvg$ is a Lie group and $\fczvg$ is the component of $\fvg$
that contains $\one$.
Let $\Omega$ be a component of $\fvg$.
If (\ref{fxvb}) holds and if there exists $x_0 \in \Omega$ such that $f(x_0)=x_0(\vb)$,
then $f(x)=x(\vb)$ for any $x \in \Omega$ \label{ptcomponenent}.


\section{Range decreasing additive maps \label{additive}}

The main result of this section is the following

\begin{thm} \label{fve}
  Let $\sce$ be a locally convex space with $\dim_{\br} \sce \ge 2$,
  and let $f: \fve \to \sce$ be a range decreasing additive group homomorphism
  such that $ f|_{\fcve}$ $\not\equiv$ $\zero$. Then $f$
  is the evaluation $E_{\vb}$ at some point $\vb$ of $V$.
\end{thm}
\begin{proof}
  By Lemma \ref{rde}, there exists $\vb \in V$ such that $f(x)=x(\vb)$,
  $x \in \fcve$.
  We say that the range of $x_0 \in \fve$ is contained in an open half-space
  of $\sce$ if there exist a nonzero continuous linear functional $\mu_0: \sce \to \br$
  and $t_0 \in \br$ such that
  $x_0(V) \subset \{a \in \sce: \mu_0(a) >t_0 \}$.
  Assume for contradiction that $f(x_0)\not=x_0(\vb)$.
  Take $\chi_0 \in C^{\infty}_c(V, \br)$ such that $\chi_0$ is constant $1$ on a neighborhood of $\vb$.
  Then
  $f(x_0-\chi_0 x_0)=f(x_0)-x_0(\vb) \not=\zero$.
  Note that $\mu_0$ is bounded on the compact subset $(x_0-\chi_0 x_0)(\supp \chi_0)$. Thus there exists
  $t'_0 \in \br$ such that
  $$(x_0-\chi_0 x_0)(V) \subset \{a \in \sce: \mu_0(a) >t'_0 \}.$$
  Without loss of generality, we may assume that $x_0$ vanishes on an open neighborhood
  of $\vb$ and $f(x_0) \not=\zero$, otherwise replace it by $x_0-\chi_0 x_0$.
  Take $ a_0 \in \sce \setminus \left(x_0(V) \cup \{\br f(x_0) \}\right)$
 and $\chi_1 \in C^{\infty}_c(V, [0, 1])$
 such that $\chi_1$ is constant $1$ on a neighborhood of $\vb$ and
 $\supp \chi_1 \subset \{v \in V: x_0(v)=\zero \}$. Note that $a_0-f(x_0) \not=\zero$. Set
 $x_1=\chi_1 \left(a_0-f(x_0) \right) \in \fcve$.
Then there exists
$v_1 \in V$ such that
$a_0=f(x_0+x_1)=x_0(v_1)+x_1(v_1)$.
Since $a_0 \not\in x_0(V)$, we have $\chi_1(v_1) \not=0$. Hence $x_0(v_1)=\zero$.
Therefore $a_0=\chi(v_1)(a_0-f(x_0)),$
which implies that $a_0 \in \{\br f(x_0) \}$. We have a contradiction.
So $f(x_0)=x_0(\vb)$.

Let $a_1 \in \sce \setminus \{\zero\}$, $\mu_1: \sce \to \br$ a nonzero continuous linear functional
and $\tau \in C^{\infty}(\br, [0, 1])$ such that $\mu_1(a_1)=1$, $\tau(t)=1$ for $t \ge 1$ and $\tau(t)=0$ for $t \le 0$.
For any $x_2 \in \fve \setminus \fcve$, we define an $\cf$ map $V \to \sce$ by
$$x_{2, 1}=\left(\tau \circ \mu_1 \circ x_2-1 \right) (\mu_1 \circ x_2) a_1+x_2.$$
Put
$x_{2, 2}=x_2-x_{2, 1} \in \fve$.
It is straightforward to verify that $\mu_1 \circ x_{2, 1}(v) \ge 0$ and
 $\mu_1 \circ x_{2, 2}(v) \le 1$ for each $v \in V$. So we conclude that $f(x_2)=f(x_{2, 1})+f(x_{2, 2})=x_2(\vb)$.
\end{proof}

The long line $\mathbb{L}$ \label{longline} is a connected
one dimensional  Hausdorff topological manifold that possesses
uncountably many $C^{\infty}$ structures
 \cite{ny}. It is not realcompact \cite[Proposition 3.3.5]{mar}, \cite[Remark 1.2]{da}.
Consequently, there exists
a nonzero algebra homomorphism
$T: C(\mathbb{L}, \br) \to \br$
that is not simply the evaluation at a specific point in $\mathbb{L}$.
Define
an algebra homomorphism $T_{\bc}: C(\mathbb{L}, \bc) \to \bc$
by
$T_{\bc}(x+iy)=T(x)+iT(y)$, where $x, y \in C(\mathbb{L}, \br)$.
It follows from Proposition \ref{algebrard} that $T_{\bc}$ is range decreasing.
By Theorem \ref{fve}, the restriction $T|_{C_c(\mathbb{L}, \br)}$
is identically zero.

Theorem \ref{fve} does not hold when the locally
convex space $\sce$ is replaced by a general connected Lie group $G$.
The Lie group $\mathrm{SU}(2)$ is isomorphic to the Lie group of unit quaternions. As
a manifold, we can identify $\mathrm{SU}(2)$ with
$S^3$. The space $C^{\infty}(S^3, \mathrm{SU}(2))$
is a Fr\'echet Lie group.
We denote the topological degree of a map $x \in C^{\infty}(S^3, \mathrm{SU}(2))$ as $\sd(x)$.
If $\sd(x) \not=0$, then $x(S^3)=\mathrm{SU}(2)$.
According to the Hopf's Degree Theorem, two maps $x_1, x_2 \in C^{\infty}(S^3, \mathrm{SU}(2))$
belong to the same component if and only if $\sd(x_1)=\sd(x_2)$.
Let $x_{id}: S^3 \to \mathrm{SU}(2)$
be the identity map.  We write $x_{id}^j \in C^{\infty}(S^3, \mathrm{SU}(2))$ for the map
$ S^3 \ni \cm \mapsto \cm^j \in \mathrm{SU}(2)$, $j \in \bz$.

 In matrix notation, the group $\mathrm{SU}(2)$ consists of the matrices
\begin{equation*}
\cm(\alpha, \beta)=\begin{pmatrix} \alpha & -\overline{\beta}\\\beta &
\overline{\alpha}
\end{pmatrix}, \hspace{2mm} \alpha, \beta \in \bc, \hspace{2mm} |\alpha|^2+|\beta|^2=1.
\end{equation*}
Let $T$ be the  maximal torus $\{\cm(e^{i\theta}, 0) \in \mathrm{SU}(2): \theta \in \br \}$.
The centralizer of
$\cm(e^{i\theta}, 0)$, where $\theta \not=k\pi$, $k \in \bz$, in $\mathrm{SU}(2)$ is exactly $T$.
Moreover, since any matrix
 $\cm \in \mathrm{SU}(2) \setminus \{\pm \one \}$
 belongs to the conjugacy class of some $\cm(e^{i\theta}, 0)$ with $\theta \not=k\pi$,
it follows that its centralizer $Z(\cm)$
 is a maximal torus. This conclusion is based on the action of
  $\mathrm{SU(2)}$ on itself by conjugation, where the
  stabilizers of distinct points within the same orbit are conjugate to one another.
Note that $Z(\cm)$ is the only maximal torus
that contains $\cm$. Additionally,  $Z(\cm)$
is an invariant subset of the maps $x_{id}^j$ for $j \in \bz$.
Thus any point in $\mathrm{SU}(2) \setminus \{\pm \one \}$
has exactly $|j|$ inverse images under the map $x_{id}^j$.

Let
$\exp_{\mathrm{SU(2)}}: \mathrm{su}(2) \to \mathrm{SU(2)}$ be the exponential map of $\mathrm{SU(2)}$.
We can find an open neighborhood $\cd$ of $\zero \in \mathrm{su(2)}$
such that $\exp_{\mathrm{SU(2)}}$ maps it
 diffeomorphically onto $\mathrm{SU(2)} \setminus \{-\one\}$ \cite[Section 4.2]{st}.
Using the local coordinates induced by $\exp_{\mathrm{SU(2)}}^{-1}$,
it is straightforward to verify that the differential of $x_{id}^j$ at $\cm \in \mathrm{SU(2)}$, where $\cm^j \not=\pm \one$,
preserves orientation when $j>0$ and reverses it
when $j<0$.
Thus we have
$\sd(x_{id}^j)=j, \hspace{1.5mm} j \in \bz$. This leads to the conclusion that
$\sd(x_1 x_2)=\sd(x_1)+\sd(x_2)$, $x_1, x_2 \in C^{\infty}(S^3, \mathrm{SU}(2))$.
We have  established the following

\begin{prop} \label{su2}
  Let $f_0$ be the group isomorphism
  $$C^{\infty}(S^3, \mathrm{SU}(2)) \ni x \mapsto  (-\one)^{\sd(x)}x \in  C^{\infty}(S^3, \mathrm{SU}(2)). $$
  For each $v \in S^3$,
  $E_v \circ f_0$
  is a range decreasing group homomorphism. However, it does not correspond to the evaluation at a specific point in
 $S^3$.
\end{prop}

\section{The non-surjective elements of $\fvg$ \label{rdcomposition}}

The main result of this section is the following

\begin{thm} \label{rdgroup}
  Let $f: \fcvg \to G$ (respectively $f: \fvg \to G$) be a range decreasing group homomorphism
  with $ f|_{\fczvg} \not\equiv \one$.
   If $\dim_{\br} G \ge 2$,
  then there exists $\vb \in V$ such that $f(x)=x(\vb)$ for any $x$ in the domain of $f$
  with $x(V) \not= G$.
\end{thm}

If $G$ is not compact, then
any map $x \in \fcvg$ cannot be surjective. Furthermore,
if $\dim_{\br} G >\dim_{\br} V$ and $V$ is second countable,
then according to the Morse-Sard theorem, any $C^1$ map $V \to G$
is also non-surjective.  In the case where $V$ is $\sigma$-compact and
$G$ is an infinite dimensional Banach Lie group,
any continuous map $x: V \to G$  will similarly fail to be surjective.
Theorem \ref{rdgroup} offers a stronger result than \cite[Theorem 1.1]{z23}, even
in scenarios where $V$ is compact.

\begin{lemma} \label{ymy1}
  Let $\vb \in V$, $x_0 \in \fvg$, $U \subset G$ an open
  neighborhood of $\one \in G$ and $a_0, a_1 \in  G \setminus \{x_0(\vb)\}$.
  If $\dim_{\br} G \ge 2$, then
  there exist a neighborhood $\tilde{V}$ of $\vb$ and finitely many
  $y_1, \cdots, y_m \in C^{\infty}_c(V, U) \subset \fvg$ such that
  $\supp y_j \subset \tilde{V}$, $j=1, \cdots, m$,
  \begin{eqnarray}
  & a_1=y_m(\vb) \cdots y_1(\vb)a_0  \hspace{2mm} \text{and} & \label{zeta1} \\
  & a_1 \not\in  y_m \cdots y_1  x_0 ( \tilde{V} ). & \label{y1y2}
  \end{eqnarray}
\end{lemma}
\begin{proof}
Take a $C^{\infty}$ path $\Gamma: [0, 1] \to G \setminus \{a_1\}$
  with
  $\Gamma(0)=x_0(\vb)$ and $\Gamma(1)=a_1 a_0^{-1} x_0(\vb)$.
Let $\exp_G: \fg \to G$ be the exponential map of $G$,
and let $\cd$ be an open neighborhood of $\zero \in \fg$
such that $\exp_G|_{\cd}: \cd \to \exp_G(\cd)$ is a diffeomorphism and $\exp_G(\cd) \subset U$.
Choose an open convex neighborhood
$\cd' \subset \cd$
of $\zero \in \fg$,
an open neighborhood $\tilde{V}$ of $\vb$ and $c_0>0$ such that
\begin{eqnarray} \label{zetanot}
& a_1 \not\in \Gamma(t) \exp_G(\cd'), \hspace{1.5mm} t \in [0, 1], \hspace{1.5mm} x_0(\tilde{V})
\subset x_0(\vb) \exp_G(\cd') \hspace{1.5mm} \text{and} & \\
 &  \Gamma(t_1) \Gamma(t_2)^{-1} \in \exp_G(\cd') \hspace{1mm} \hspace{1.5mm} \text{when} \hspace{1.5mm} |t_1-t_2|<c_0. & \label{c0}
   \end{eqnarray}
   Let $\kappa \in C^{\infty}_c(V, [0, 1])$ be a function such that $\kappa(\vb)=1$ and $\supp \kappa \subset \tilde{V}$.
   Additionally, consider
    a partition $0=\tilde{t}_0<\tilde{t}_1<\cdots<\tilde{t}_m=1$ of the interval $[0, 1]$
   with $\tilde{t}_j-\tilde{t}_{j-1}<c_0$, $j=1, \cdots, m$.
     Define $y_j \in C^{\infty}_c(V, G)$ by
  $$y_j(v)=\Gamma\left(\kappa(v)\tilde{t}_j \right)
  \Gamma\left(\kappa(v)\tilde{t}_{j-1} \right)^{-1} \in G, \hspace{1mm} v \in V, \hspace{1mm} j=1, \cdots, m.$$
  It is clear that $\supp y_j \subset \tilde{V}$.
  By (\ref{c0}), we have $y_j(V) \subset \exp_G(\cd') \subset U$.
  It is straightforward to verify that (\ref{zeta1}) holds.
   Note that
   $y_m(v) \cdots y_1 (v) x_0(\vb) \in \Gamma([0, 1])$, $v \in V$.
   Application of (\ref{zetanot}) gives (\ref{y1y2}).
\end{proof}

\noindent {\it Proof of Theorem \ref{rdgroup}. \label{pfrdgroup}} By Theorem \ref{nbhd1}, there exists $\vb \in V$ such that
$f(\exp_G \circ \xh)=\exp_G \circ \xh(\vb)$, $\xh \in \cf_c(V, \fg)$.
In particular, there exists an open neighborhood $U$ of $\one \in G$ such that $f(x)=x(\vb)$ for any
  $ x \in \cf_c(V, U)$.
Suppose that there exists $x_0 \in \fcvg$ (respectively  $x_0 \in \fvg$) such that $x_0(V) \not=G$ and $f(x_0) \not=x_0(\vb)$.
By Lemma \ref{ymy1}, where we choose $a_0=f(x_0)$ and $a_1 \in G \setminus x_0(V)$,
we can find a neighborhood $\tilde{V}$ of $\vb$ and finitely many $y_1, \cdots, y_m \in C^{\infty}_c(V, U)$ such that
  $\supp y_j \subset \tilde{V}$, $j=1, \cdots, m$,
  $$a_1=y_m(\vb) \cdots y_1(\vb) f(x_0)=f(y_m \cdots y_1 x_0)$$
  and $a_1 \not\in  y_m \cdots y_1  x_0 ( \tilde{V} )$.
  For any $v \in V \setminus \tilde{V}$,
  $ y_m(v) \cdots y_1(v)  x_0(v)=x_0(v) \not=a_1.$
  Hence $a_1 \not\in y_m \cdots y_1  x_0(V)$,
  which contradicts the range decreasing condition of $f$. \qed

\begin{cor} \label{s1}
Let $S$ be a nonempty set, $\sg$ the group of all maps $S \to G$
  and $f: \fczvg \to \sg$ a nonconstant range decreasing group homomorphism.
  Put $S_f=\{s \in S: E_s \circ f \not\equiv \one \}$.
    If $\dim_{\br} G=1$, we additionally require that for any $s \in S_f$,
    the set $\cs_{E_s \circ f}$ consists of a single point.
  Then there exists a map $\phi: S_f \to V$ such that
  \begin{equation} \label{fxs}
    f(x)(s)=
    x\left(\phi(s) \right), \hspace{2mm} s \in S_f,
  \end{equation}
  for any $x \in \fczvg$. Furthermore,
  if $S$ is a finite dimensional $C^{\infty}$ manifold $W$,
  possibly with boundary, and $f(x) \in \tilde{\cf}(W, G)$ for any $x \in \fczvg$, then
  $S_f$ is an open subset of $W$ and $\phi$ is an $\tilde{\cf}$ map.
\end{cor}
\begin{proof}
  Application of Theorem \ref{nbhd1} to $E_s \circ f$, where $s \in S_f$, gives (\ref{fxs}).
  If $f(x) \in \tilde{\cf}(W, G)$ for any $x \in \fczvg$, it
  follows from the continuity of $f(x)$ that
  $S_f$ is open.
  By (\ref{fxs}) and Proposition
  \ref{phi}, $\phi$ is an $\tilde{\cf}$ map.
\end{proof}

\begin{cor} \label{rdgh}
  Let $f: \fcvg \to \sg$ (respectively $ f: \fvg \to \sg$) be a range decreasing group homomorphism
  with $f|_{\fczvg} \not\equiv \one$. Put $S_f=\{s \in S: E_s \circ f|_{\fczvg} \not\equiv \one \}$.
  If $\dim_{\br} G \ge 2$, then (\ref{fxs}) holds
  for any $x$ in the domain of $f$ with $x(V) \not=G$. Furthermore,
   if $G$ is a locally convex space $\sce$ with $\dim_{\br} \sce \ge 2$,
  then (\ref{fxs}) holds for any $x$ in the domain of $f$.
\end{cor}
\begin{proof}
  Apply Theorems \ref{rdgroup}, \ref{fve} and Lemma \ref{rde}
  to $E_s \circ f$, $s \in S_f$.
\end{proof}


\section{Automatic continuity and automatic smoothness \label{wco1}}

The main result of this section is the following

\begin{thm} \label{frechetlg}
 Let $\pi_{\scg_1, V}: \scg_1 \to V$ (respectively  $\pi_{\scg_2, W}: \scg_2 \to W$) be a $C^{\infty}$ Banach Lie group
bundle with the typical fiber $G_1$ (respectively $G_2$), and
 let $f: \fczvscgo \to \ftscgt$ be a group homomorphism.
 If there exist
 a map $\phi: W \to V$ and a (not necessarily continuous) bundle morphism
   $\gamma: \phi^{\ast} \scg_1 \to \scg_2$ over the identity map on $W$
   such that $\phi$ is not constant on any open subset of $W$, the restriction of $\gamma$ to each fiber
   is a nonconstant group homomorphism and
    \begin{equation} \label{gwco}
       f(x)=\gamma \circ (\phi^{\ast} x), \hspace{2mm} x \in \fczvscgo,
    \end{equation}
 then $\gamma$  is continuous
 and $\phi$ is an $\tilde{\cf}$ map.
\end{thm}

If the map $\phi$ in Theorem \ref{frechetlg} is constant, \label{disc}
the bundle morphism $\gamma$ may exhibit discontinuity. For instance,
let $v_0 \in V$, and let
$\gamma_0$ be a discontinuous group homomorphism $G_1 \to G_2$. We can then
examine the group homomorphism $\cf_c^0(V, G_1) \ni x \mapsto \gamma_0(x(v_0)) \in G_2 \subset \tilde{\cf}(W, G_2)$.

\begin{lemma} \label{polish}
   Let $G$ be a metrizable topological group
  and
  $\ch$ a locally compact Polish topological group. Then any Borel measurable
  homomorphism $\psi: \ch \to G$ is continuous. In particular, if $\psi_0: \ch \to G$ is a map
  such that there exists a sequence of continuous homomorphisms $\{\psi_j: \ch \to G\}$ with
  $\psi_0(a)=\lim_{j \to \infty} \psi_j(a)$, $a \in \ch,$
  then $\psi_0$ is also a continuous homomorphism.
\end{lemma}

If $G$ is either separable or $\sigma$-compact, then Lemma \ref{polish} is a consequence of
\cite[Theorem 22.18]{hewitt}. Additionally, if $G$ is locally compact, then Lemma \ref{polish}
follows from \cite[Theorem 1]{kl}.

\vspace{2mm}

\noindent {\it Proof of Lemma \ref{polish}. }
Let $\lambda$ be the left
  Harr measure on $\ch$, which is a Radon measure. Note that $\ch$ is $\sigma$-compact.
  Thus $\lambda$ is locally determined \cite[211L]{fr2}, and
  $\ch$ is a Radon measure space \cite[P. 388]{fr} (where we may
  replace $\lambda$ by its completion).
  It follows from \cite[Theorem 2B]{fr} that $\psi$ is almost continuous. Any Radon measure
  on a $\sigma$-compact space is regular.
  So there exists a compact subset $K \subset \ch$ such that $\lambda(K)>0$ and
  $\psi$ is continuous on $K$.

  Note that $K K^{-1}$ is a neighborhood of $\one \in \ch$ \cite{s}.
  We claim that $\psi$ is continuous on $K K^{-1}$, which implies that $\psi$ is continuous on $\ch$.
  Otherwise there exist $a_0 b_0^{-1} \in K K^{-1}$ and a sequence $\{a_j b_j^{-1}\} \subset K K^{-1}$
  such that $\lim_{j \to \infty} a_j b_j^{-1} \to a_0 b_0^{-1}$
  and $\lim_{j \to \infty} \psi(a_j b_j^{-1}) \not=\psi(a_0 b_0^{-1})$. Take
  an open neighborhood $O$ of $\psi(a_0 b_0^{-1})$ and a subsequence
  $\{a_{j_k} b_{j_k}^{-1}\}$ of $\{a_j b_j^{-1}\}$ such that $\psi(a_{j_k} b_{j_k}^{-1}) \not\in O$, $k=1, 2, \cdots$.
  Since $K$ is compact, we may assume that the sequences $\{a_{j_k}\}$ and $\{b_{j_k}\}$ in $K$
  are convergent (otherwise replace them by convergent subsequences). Suppose that
  $a_{j_k} \to \alpha_0 \in K$ and $b_{j_k} \to \beta_0 \in K$ when $k \to \infty$. Then $\alpha_0 \beta_0^{-1}=a_0 b_0^{-1}$
  and
  $\lim_{k \to \infty} \psi(a_{j_k} b_{j_k}^{-1})=\psi(\alpha_0) \psi(\beta_0)^{-1}=\psi(a_0 b_0^{-1})$.
  We have a contradiction.

  By Proposition \ref{measurable}, the map $\psi_0$ in  the lemma is a Borel measurable group homomorphism.
  Hence $\psi_0$ is continuous.
\qed

\begin{lemma} \label{passsequence}
  Let $G$ be a connected locally exponential  Fr\'echet Lie group, $v_0 \in V$ and
  $\{v_k\} \subset V \setminus \{v_0\}$ a pairwise distinct sequence with $\lim_{k \to \infty} v_k=v_0$.
  Then there exists a sequence of open neighborhoods $\{U_k\}$ of $\one \in G$ such that for any sequence
  $\{a_k\} \subset G \setminus \{\one\}$ with $a_k \in U_k$, $k \in \bn$,
  we have $\lim_{k \to \infty} a_k =\one$,
  and we can find a $C^{\infty}$ map $x_0 \in \fczvg$ that satisfies the conditions
   $x_0(V) \not=G$ and $x_0(v_k)=a_k$, $k \in \bn$
   (hence $x_0(v_0)=\one$).
\end{lemma}
\begin{proof} Let $\exp_G: \fg \to G$ be the exponential map of $G$
and $\cd$ an open convex neighborhood of $\zero \in \fg$
such that the map $\exp_G|_{\cd}: \cd \to \exp_G(\cd)$ is a diffeomorphism.
Consider a precompact open neighborhood $\tv$ of $v_0$ such that
there exists a local coordinate chart $(U, \Phi)$ for $V$
with $\overline{\tv} \subset U$. Take $N \in \bn$ with $v_k \in \tv$, $k>N$.
   Choose cut-off functions
   $\chi_k \in C^{\infty}_c(\tv, [0, 1])$, $k >N$,
   such that $\Phi(\supp \chi_k)$
   is contained in an open ball $B_k \subset \Phi(\tv)$ centered at $\Phi(v_k)$,
   $\chi_k(v_k)=1$,
  $\Phi(v_0) \not\in \overline{B_k}$
   and $\overline{B_k} \cap \overline{B_j}=\emptyset$ when $k \not=j$, $k, j >N$.
   Let $\mathscr{P}=\{p_1, p_2, \cdots\}$ be a family of seminorms that defines the topology on $\fg$.
   Set
   $\hat{U}_k=\{\ah \in \fg: p_i(\ah)<r_k, i=1, 2, \cdots, k \}$, $k>N$, where
   $r_k =2^{-k}||\chi_k \circ \Phi^{-1}||_{C^k}^{-1}$,
   $||\chi_k \circ \Phi^{-1}||_{C^k}=\max _{|\alpha|\le k} \max_{t \in \Phi(\overline{\tv})} \left|\partial^{\alpha} (\chi_k  \circ \Phi^{-1})(t) \right|$,
   $t=(t_1, \cdots, t_{n}) \in \br^{n}$, $n=\dim_{\br} V$, $\alpha=\left(\alpha_1, \cdots, \alpha_n\right)$,
   $\alpha_1, \cdots, \alpha_n \in \bn \cup \{0\}$, $|\alpha|=\alpha_1+\cdots+\alpha_n$ and
   $$ \te
    \partial^{\alpha} (\chi_k  \circ \Phi^{-1})(t)=\frac{\partial^{|\alpha|}}{\partial t_1^{\alpha_1}  \cdots \partial t_n^{\alpha_n}} (\chi_k  \circ \Phi^{-1})(t).
   $$
   Note that
   $||\chi_k \circ \Phi^{-1}||_{C^k} \ge 1$.
   By possibly increasing $N$,  we can assume that
   $\hat{U}_k \subset \frac{1}{2}\cd$, $k>N$.
   Let
   \begin{equation} \label{uk}
   U_k=\exp_G(\hat{U}_k), \hspace{2mm} k>N, \hspace{2mm} U_1=\cdots=U_N=\exp_G(\cd/4),
   \end{equation}
 and let $\{a_k\} \subset G$ be a sequence with $a_k \in U_k$.
 Set $\ah_k=(\exp_G|_{\cd})^{-1}(a_k)$, $k \in \bn$.
 Since $\lim_{k \to \infty} p_i(\ah_k)=0$ for any seminorm $p_i \in \mathscr{P}$,
   it follows that $\lim_{k \to \infty} \ah_k=\zero$. Define a map
   $$\te \tx_0: V \ni v \mapsto \sum_{k=N+1}^{\infty} \chi_k(v) \ah_{k} \in \fg.$$
    For any $v \in V$, the series contains at most one term that does not vanish at $v$.
    Note that $\tx_0(V) \subset \frac{1}{2} \cd$ and
   $\tilde{x}_0(v_k)=\ah_{k}$, $k>N$.
   For any continuous linear functional $\Lambda: \fg \to \br$,
   we can find a constant $c>0$ and a finite collection of seminorms $p_{i_1}, \cdots, p_{i_m} \in \mathscr{P}$
   such that
   $$|\Lambda(\ah)| \le c \max\{p_{i_1}(\ah), \cdots, p_{i_m}(\ah) \}, \hspace{2mm} \ah \in \fg.$$
   It is clear that
   $\Lambda \circ \tilde{x}_0(v)=\sum_{k=N+1}^{\infty} \chi_k(v) \Lambda(\ah_{k})$
   is a $C^{\infty}$ function $V \to \br$. Thus $\tilde{x}_0 \in C^{\infty}_c(V, \fg)$
   \cite{fro}.
   Finally, we construct the map $x_0$ by modifying the value
   of $\exp_G \circ \tilde{x}_0$ in small neighborhoods of the points $v_1, \cdots, v_N$.
   \end{proof}

\begin{lemma} \label{phifp}
  Let $G_1$, $G_2$ be connected locally exponential
  locally convex Lie groups and $f: \fczvgo$ $\to$ $\tilde{\cf}(W, G_2)$ a group homomorphism.
  If there exist
   maps $\phi: W \to V$ and
  $\gamma: W \times G_1 \to G_2$
  such that $\gamma(w, \cdot)$
  is a nonconstant group homomorphism for any $w \in W$ and
  \begin{equation} \label{fcz}
    f(x)(w)=\gamma(w, x \circ \phi(w)), \hspace{2mm} w \in W, \hspace{2mm} x \in \fczvgo, \hspace{2mm} \text{then}
  \end{equation}
  \begin{equation} \label{fp}
  \gamma(\cdot, a) \in \tilde{\cf}(W, G_2), \hspace{2mm} a \in G_1.
  \end{equation}
  If we further assume that $\gamma$ is continuous,
   then $\phi$ is an $\tilde{\cf}$ map.
\end{lemma}
\begin{proof}
  By Proposition \ref{phicontinuous}, $\phi$ is continuous.
  Let $\exp_{G_i}: \fg_i \to G_i$ be the exponential map of $G_i$, $i=1, 2$,
    and $\cd_i$ an open convex neighborhood of $\zero \in \fg_i$
    such that $\exp_{G_i}|_{\cd_i}: \cd_i \to \exp_{G_i}(\cd_i)$ is a diffeomorphism,
  $\ah \in \fg_1 \setminus \{\zero\}$ and $w_0 \in W$. Choose a function
  $\chi_{\phi(w_0)} \in C^{\infty}_c(V, [0, 1])$ such that  $\chi_{\phi(w_0)}$ is constant
  $1$ on a neighborhood of $\phi(w_0)$. Set
  $x_{\ah, \phi(w_0)}=\exp_{G_1}(\chi_{\phi(w_0)} \ah) \in \fczvgo$.
  Then there exists a
  neighborhood $O_{w_0}$ of $w_0$ such that
  \begin{equation} \label{gamma}
  \gamma(w, \exp_{G_1}(\ah))=\gamma(w, x_{\ah, \phi(w_0)} \circ \phi(w) )=f(x_{\ah, \phi(w_0)})(w), \hspace{2mm} w \in O_{w_0}.
  \end{equation}
  The right hand side of the equation above is an $\tilde{\cf}$ map $O_{w_0} \to G_2$.
  Since $w_0$ is arbitrary, we have (\ref{fp}) for $a \in \exp_{G_1}(\fg_1)$.
  Since $G_1$ is generated by the subset $\exp_{G_1}(\fg_1)$ \cite[Theorem 7.4]{hewitt},
  (\ref{fp}) is valid for any $a \in G_1$.

  If $\gamma$ is continuous, then for any $w_0 \in W$,
  there exist an open neighborhood $W_0$ of $w_0$ and an open convex neighborhood $U \subset \cd_1$
  of $\zero \in \fg_1$ such that
  $\gamma(w, a) \in \exp_{G_2}(\cd_2)$ for every $(w, a) \in W_0 \times \exp_{G_1}(U)$.
  Let $\tilde{\gamma}(w, \cdot): \fg_1 \to \fg_2$ be the Lie algebra homomorphism of $\gamma(w, \cdot)$,
  where $w \in W$.
  Take $\ah_0 \in U  \setminus \{\zero\}$ and a continuous linear functional $\mu: \fg_2 \to \br$
  such that
  $\mu \circ \tilde{\gamma}(w_0, \ah_0) \not=0$.
  Note that $\tilde{\gamma}(w, \ah_0)=(\exp_{G_2}|_{\cd_2})^{-1} \circ \gamma(w, \exp_{G_1}(\ah_0))$, where $w \in W_0$.
   According to (\ref{fp}), $\tilde{\gamma}(\cdot, \ah_0) \in \tilde{\cf}(W, \fg_2)$.
  By shrinking $W_0$ if necessary,
  we may assume that
  $ \mu \circ \tilde{\gamma}(w, \ah_0) \not=0$ for each $w \in W_0$.
  For any
  $\chi \in C^{\infty}_c(V, [0, 1])$, we have
  $ \exp_{G_1}(\chi \ah_0) \in \fczvgo$.
  It follows from (\ref{fcz}) that
  $$\chi \circ \phi(w) \tilde{\gamma}(w, \ah_0) =\tilde{\gamma}(w, \chi \circ \phi(w) \ah_0)
  = (\exp_{G_2}|_{\cd_2})^{-1} \left( f(\exp_{G_1}(\chi \ah_0))(w) \right),$$
  where $ w \in W_0$. Therefore
  $$\chi \circ \phi(w)=
  \mu \circ (\exp_{G_2}|_{\cd_2})^{-1} \left( f(\exp_{G_1}(\chi \ah_0))(w) \right)/\mu \circ \tilde{\gamma}(w, \ah_0) \in \br$$
  is an $\tilde{\cf}$ function on $W_0$, which implies that $\phi$ itself is an $\tilde{\cf}$ map on $W_0$.
  \end{proof}

\begin{thm} \label{wco2}
  Let $G_1$ and $G_2$ be connected Banach Lie groups, and let
  $f: \fczvgo$ $\to$ $\tilde{\cf}(W, G_2)$ be a group homomorphism.
  If there exist
   maps $\phi: W \to V$ and
  $\gamma: W \times G_1 \to G_2$
  such that $\phi$ is not constant on any open subset of $W$,
  $\gamma(w, \cdot)$
  is a nonconstant group homomorphism for any $w \in W$ and
    $f(x)(w)=\gamma(w, x \circ \phi(w))$, $w \in W$, $x \in \fczvgo$,
      then
      $\gamma(w, \cdot)$ is a Lie group homomorphism for any $w \in W$.
  Let $\tilde{\gamma}(w, \cdot): \fg_1 \to \fg_2$ be the Lie algebra homomorphism of
  $\gamma(w, \cdot)$. Then $\tilde{\gamma}(\cdot, \ah): W \to \fg_2$
  is an $\tilde{\cf}$ map for any $\ah \in \fg_1$.
\end{thm}
\begin{proof}
 Let $\exp_{G_i}: \fg_i \to G_i$ be the exponential map of $G_i$,  where $i=1, 2$,
    and let $\cd_i$ be an open convex neighborhood of $\zero \in \fg_i$
    such that $\exp_{G_i}|_{\cd_i}: \cd_i \to \exp_{G_i}(\cd_i)$ is a diffeomorphism.
    Note that $\fg_i$ is a Banach space.
Let $w_0 \in W$, and consider a sequence $\{w_k\}$ in $W \setminus \{w_0\}$
such that
   $\lim_{k \to \infty} w_k$ $=$ $w_0$ and the points
   $\phi(w_0)$, $\phi(w_1)$, $\phi(w_2)$, $\cdots $
   are all distinct from one another.

Assume for contradiction that $\gamma(w_{k}, \cdot)$ is discontinuous for each $k \in \bn$.
   Let $B \subset \frac{1}{2}\cd_2$
   be an open ball centered at $\zero \in \fg_2$, and let the neighborhoods $U_k$, $k \in \bn$,
   of $\one \in G_1$ be as in Lemma \ref{passsequence}.
   We may assume that $U_k= \exp_{G_1}(\hat{U}_k)$, $k \in \bn$,
   where $\hat{U}_k \subset \cd_1$ is an open convex neighborhood of $\zero \in \fg_1$
   ((\ref{uk})).
   We can find an open ball $B_k \subset B$
   centered at $\zero \in \fg_2$ and a sequence $\{\ah_{k, j} \} \subset \hat{U}_k \setminus \{\zero\}$
   such that
   $$\te \lim_{j \to \infty} \ah_{k, j}=\zero \hspace{2mm} \text{and} \hspace{2mm}  \gamma(w_{k}, \exp_{G_1}(\ah_{k, j})) \not\in \exp_{G_2}(B_k), \hspace{2mm} j, k \in \bn.$$
   For a given $k \in \bn$, if $\gamma(w_{k}, \exp_{G_1}(\ah_{k, j_1})) \not\in \exp_{G_2}(B)$ for some $j_1 \in \bn$,
   set
   $$ a_k=\exp_{G_1}(\ah_{k, j_1}) \in U_k \setminus \{\one\}.$$
   Otherwise $\gamma(w_{k}, \exp_{G_1}(\ah_{k, j})) \in \exp_{G_2}(B \setminus B_k)$ for each $j \in \bn$.
    Let $m_k \in \bn$ be chosen such that $m_k \bhat \not\in B$ for every
   $\bhat \in B \setminus B_k$.  Additionally, select $\ah_{k, j_2} $ from the sequence $\{\ah_{k, j} \}$
   such that
   $ m \ah_{k, j_2} \in \hat{U}_k \setminus \{\zero\}$, $m=1, \cdots, m_k$.
   Let $\bhat_{k, j_2}$ be the unique element in $B \setminus B_k$ with
   $\exp_{G_2}(\bhat_{k, j_2})=\gamma(w_{k}, \exp_{G_1}(\ah_{k, j_2}))$.
   Then there is a positive integer $n_k \le m_k$ such that $(n_k-1) \bhat_{k, j_2} \in B$ and $n_k \bhat_{k, j_2} \in \cd_2 \setminus B$.
   In this case, set
   $$a_k=\exp_{G_1}(n_k\ah_{k, j_2}) \in U_k \setminus \{\one\}.$$
    If $\gamma(w_k, a_k) \not\in \exp_{G_2}(\cd_2)$, then $\gamma(w_k, a_k) \not\in \exp_{G_2}(B)$.
   If $\gamma(w_k, a_k) \in \exp_{G_2}(\cd_2)$, we still have
   $$\gamma(w_k, a_k)=\gamma(w_{k}, \exp_{G_1}(\ah_{k, j_2}))^{n_k}=\exp_{G_2}(n_k\bhat_{k, j_2}) \not\in \exp_{G_2}(B).$$
   Therefore we obtain a sequence $\{a_k\} \subset G_1 \setminus \{\one\}$ with $a_k \in U_k$
   and $\gamma(w_k, a_k) \not\in \exp_{G_2}(B)$ for every $k \in \bn$. By Lemma \ref{passsequence},
   there exists $x_0 \in \fczvgo$ such that $x_0 \circ \phi(w_k)=a_k$ for each $k \in \bn$
   and $x_0 \circ \phi(w_0) =\one$.
  Thus
  $$\gamma(w_k, a_k)=f(x_0)(w_k) \to f(x_0)(w_0)=\gamma(w_0, x_0 \circ \phi(w_0))=\one$$
  when $k \to \infty$, leading to a contradiction.
  Thus all but finitely many maps $\gamma(w_{k}, \cdot)$, $k \in \bn$, are continuous.
  If needed, we can adjust only a finite number of terms in the sequence
   $\{w_k\}$  to ensure that all maps
  $\gamma(w_{k}, \cdot)$, $k \in \bn$, are continuous.

   Let $\ah \in \fg_1 \setminus \{\zero\}$.
  By  the
Birkhoff-Kakutani theorem,
$G_2$ is metrizable.
  It follows from (\ref{fp}) and Lemma \ref{polish} that the group homomorphism
  \begin{equation} \label{limit}
  \te \{\br \ah\} \ni t\ah \mapsto \gamma(w_0, \exp_{G_1}(t\ah))=\lim_{k \to \infty}  \gamma(w_k, \exp_{G_1}(t\ah)) \in G_2
  \end{equation}
  is continuous.
  Let $\tilde{\gamma}(w_0, \cdot): \fg_1 \to \fg_2$ be the map whose restriction to the subspace
  $\{\br \ah\}$ is the Lie algebra homomorphism
  of the Lie group homomorphism in (\ref{limit}).

   Suppose that $\{w'_k \}$ is a sequence in $W$ such that
   $\lim_{k \to \infty} w'_k$ $=$ $w_0$ and the map $\gamma(w'_k, \cdot)$ is continuous
   for each $k \in \bn$. If the linear operators
   $\tilde{\gamma}(w'_k, \cdot)|_{\{\br \ah\}}$ for $k \in \bn$ are uniformly bounded, meaning
   there exists a constant $c>0$ such that the norm
   $||\tilde{\gamma}(w'_k, \cdot)|_{\{\br \ah\}}|| \le c$
  for each $k \in \bn$, then we can find $\delta>0$ such that $\tilde{\gamma}(w_0, t \ah) \subset \cd_2$ and
  $\tilde{\gamma}(w'_k, t \ah) \subset \cd_2$ for any $k \in \bn$ and for any $t \in (-\delta, \delta)$.
  It follows from (\ref{fp}) that
 $$\te \lim_{k \to \infty} \tilde{\gamma}(w'_k, t \ah)=\lim_{k \to \infty}
 (\exp_{G_2}|_{\cd_2})^{-1} \circ \gamma(w'_k, \exp_{G_1}(t\ah))=\tilde{\gamma}(w_0, t \ah),$$
 where $t \in (-\delta, \delta)$. Thus $\lim_{k \to \infty} \tilde{\gamma}(w'_k, \ah)=\tilde{\gamma}(w_0, \ah)$.

 Given $\ah \in \fg_1 \setminus \{\zero\}$, we will show that the linear operators
   $\tilde{\gamma}(w_k, \cdot)|_{\{\br \ah\}}$ for $k \in \bn$ are uniformly bounded.
   Assuming the contrary, we can find a subsequence $\{w_{k_i} \}$ of $\{w_k\}$
  along with a sequence $\{t_{i}\}$ in $\br$ such that $\exp_{G_1}(t_{i} \ah) \subset U_{i}$
  and $\tilde{\gamma}(w_{k_i}, t_i \ah) \subset \cd_2 \setminus B$ for each $i \in \bn$.
  According to Lemma \ref{passsequence},
   there exists a map $x_1 \in \fczvgo$ such that $x_1 \circ \phi(w_{k_i})=\exp_{G_1}(t_{i} \ah)$ for all $i \in \bn$
   and $x_1 \circ \phi(w_0) =\one$.
  Thus
  $$\gamma(w_{k_i}, \exp_{G_1}(t_{i} \ah))=f(x_1)(w_{k_i}) \to f(x_1)(w_0)=\gamma(w_0, x_1 \circ \phi(w_0))=\one$$
  as $k \to \infty$. Note that the left hand side of the equation above is contained in $\exp_{G_2}(\cd_2 \setminus B)$,
  leading to a contradiction.
 Thus we conclude that
  $ \tilde{\gamma}(w_0, \ah)=\lim_{k \to \infty} \tilde{\gamma}(w_k, \ah)$ for any $\ah \in \fg_1$.
   This implies that
   $\tilde{\gamma}(w_0, \cdot): \fg_1 \to \fg_2$ is continuous \cite[Theorem 2.8]{ru}.
   Consequently, we find that
  $\gamma(w_0, \cdot)$ is continuous for any $w_0 \in W$.

  We claim that the map $\tilde{\gamma}(\cdot, \ah): W \to \fg_2$
  is continuous. To show this, consider a sequence $\{\tilde{w}_k \}$ in $W$ with $\lim_{k \to \infty} \tilde{w}_k=w_0$.
  We can find another sequence $\{\tilde{w}'_k \}$ in $W \setminus \{w_0\}$ such that $\lim_{k \to \infty} \tilde{w}'_k=w_0$,
  $||\tilde{\gamma}(\tilde{w}_k, \cdot)|_{\{\br \ah\}}-\tilde{\gamma}(\tilde{w}'_k, \cdot)|_{\{\br \ah\}}||<2^{-k}$
  for each $k \in \bn$ and the points
   $\phi(w_0)$, $\phi(\tilde{w}'_1)$, $\phi(\tilde{w}'_2)$, $\cdots $
   are pairwise distinct. The linear operators
   $\tilde{\gamma}(\tilde{w}'_k, \cdot)|_{\{\br \ah\}}$ for $k \in \bn$ are uniformly bounded.
   As a consequence, the linear operators
   $\tilde{\gamma}(\tilde{w}_k, \cdot)|_{\{\br \ah\}}$ for $k \in \bn$ are also uniformly bounded.
   So $\lim_{k \to \infty} \tilde{\gamma}(\tilde{w}_k, \ah)=\tilde{\gamma}(w_0, \ah)$,
   which establishes the continuity of $\tilde{\gamma}(\cdot, \ah)$ at $w_0$.

  Let $O$ be a precompact open
  subset of $W$. For any fixed $\ah \in \fg_1$, there exists $j_0 \in \bn$
  such that $\tilde{\gamma}(w, \ah)/j_0 \in \cd_2$ for each $w \in \overline{O}$. Thus
  $$\tilde{\gamma}(w, \ah)=j_0 (\exp_{G_2}|_{\cd_2})^{-1} \circ \gamma(w, \exp_{G_1}(\ah/j_0)), \hspace{2mm} w \in \overline{O}.$$
   By (\ref{fp}), $\tilde{\gamma}(\cdot, \ah)$ is an $\tilde{\cf}$ map on $O$.
   \end{proof}

\noindent {\it Proof of Theorem \ref{frechetlg}.} \label{proof6.1}
The continuity of $\phi$ follows from Proposition \ref{phicontinuous}.
For any $w \in W$, we can choose an open neighborhood $O_{w} \subset W$ of $w$ and
an open neighborhood $U_{\phi(w)} \subset V$ of $\phi(w)$
such that $\phi(O_{w}) \subset U_{\phi(w)}$,
the bundle $\pi_{\scg_1, V}$ is trivial over $U_{\phi(w)}$
and the bundle $\pi_{\scg_2, W}$ is trivial over $O_{w}$.
The map  $f$ induces a group homomorphism
$$f_w: \cf_c^0(\pi_{\scg_1, V}|_{U_{\phi(w)}}) \ni x \mapsto f(x)|_{O_{w}} \in \tilde{\cf}(\pi_{\scg_2, W}|_{O_{w}}).$$
To prove the theorem, we can assume that the bundles $\pi_{\scg_1, V}$ and $\pi_{\scg_2, W}$ are trivial.
Otherwise we will consider the maps $f_w$, $w \in W$, rather than
$f$ itself.  In this case, we can then consider
$\gamma$ as a map $W \times G_1 \to G_2$.

     Next we will show that $\gamma$ is continuous.
     Hence by Lemma \ref{phifp}, $\phi$ is an $\tilde{\cf}$ map.
     Let $\exp_{G_i}: \fg_i \to G_i$ be the exponential map of $G_i$, $i=1, 2$,
    and $\cd_i$ an open convex neighborhood of $\zero \in \fg_i$
    such that $\exp_{G_i}|_{\cd_i}: \cd_i \to \exp_{G_i}(\cd_i)$ is a diffeomorphism.
     Let $w_0 \in W$, $a_0 \in G_1$ and $b_0=\gamma(w_0, a_0) \in G_2$.
     For any open neighborhood $\Omega_{b_0}$ of $b_0$,
     there exist another open neighborhood $\Omega^{\ast}_{b_0}$ of $b_0$ and
     an open neighborhood $\Omega_{\one} \subset \exp_{G_2}(\cd_2)$ of $\one \in G_2$ such that
     \begin{equation} \label{b2b1}
     b_2 b_1 \in \Omega_{b_0},  \hspace{2mm}
     b_2 \in \Omega_{\one}, \hspace{2mm} b_1 \in \Omega^{\ast}_{b_0}.
     \end{equation}
     By Theorem \ref{wco2}, the map
     $\gamma(w, \cdot)$ is a Lie group homomorphism for any $w \in W$.
     Let $\tilde{\gamma}(w, \cdot): \fg_1 \to \fg_2$
  be the Lie algebra homomorphism of $\gamma(w,\cdot)$, $\ah \in \fg_1$
  and $O_{w_0} \subset W$ a precompact open neighborhood of $w_0$.
  Theorem \ref{wco2} indicates that the subset
      $ \{\tilde{\gamma}(w, \ah): w \in \overline{O_{w_0}}\}$ of $\fg_2$
     is compact. According to the Banach-Steinhaus theorem,
     the collection of linear operators
          $\{\tilde{\gamma}(w, \cdot):  w \in \overline{O_{w_0}}\}$
      is equicontinuous.
        Thus there exists an open neighborhood $U_{\zero} \subset \cd_1$
     of $\zero \in \fg_1$ such that
     $$\gamma(w, a) \in \Omega_{\one},  \hspace{2mm} w \in \overline{O_{w_0}}, \hspace{2mm}
       a \in \exp_{G_1}(U_{\zero}).$$
     Since $\gamma(\cdot, a_0)$ is continuous ((\ref{fp})),
     there exists an open neighborhood $O'_{w_0} \subset O_{w_0}$ of $w_0$ such that
     $\gamma(w, a_0) \in \Omega^{\ast}_{b_0}$ for $w \in O'_{w_0}$.
     By (\ref{b2b1}),
     $\gamma(w, aa_0)=\gamma(w, a) \gamma(w, a_0) \in \Omega_{b_0}$ for each $w \in O'_{w_0}$
     and for each $a \in \exp_{G_1}(U_{\zero})$.
     Therefore $\gamma$ is continuous at the point $(w_0, a_0) \in W \times G_1$. \qed

     \vspace{2mm}

We  define a real vector bundle $\pi_{\eta, W}: \eta \to W$ of rank $r \in \bn$
   to be an $\cf$ bundle if its transition functions
   are $\cf$ maps from open subsets of $W$ to the general linear group $GL(r, \br)$.
   Consider two real $\cf$ vector bundles $\pi_{\eta_i, W}: \eta_i \to W$ of ranks $r_i \in \bn$ for $i=1, 2$.
   A vector bundle morphism $\kappa: \eta_1 \to \eta_2$
over the identity map on $W$
is called an $\cf$ morphism if, under local trivializations of
the bundles $\pi_{\eta_1, W}$ and $\pi_{\eta_2, W}$, $\kappa$
is represented by $\cf$ maps from open subsets of $W$ to the space of  $r_2 \times r_1$ real matrices.

\begin{prop} \label{finiterank}
    Let $\pi_{\xi_1, V}: \xi_1 \to V$ and $\pi_{\xi_2, W}: \xi_2 \to W$
   be real $C^{\infty}$ vector bundles of finite rank,
    and let
    $f: \cf_c(\pi_{\xi_1, V}) \to \tilde{\cf}(\pi_{\xi_2, W})$ be a linear map.
 If there exist
 a map $\phi: W \to V$ and a (not necessarily continuous) vector bundle morphism
   $\gamma: \phi^{\ast} \xi_1 \to \xi_2$ over the identity map on $W$
   such that the restriction of $\gamma$ to each fiber
   is nonconstant and
       $f(x)=\gamma \circ (\phi^{\ast} x)$ for each $x \in \cf_c(\pi_{\xi_1, V})$,
   then $\phi$ is an $\tilde{\cf}$ map (which implies that $\phi^{\ast} \xi_1$
   is an $\tilde{\cf}$ vector bundle)
   and $\gamma$ is an $\tilde{\cf}$ morphism.
\end{prop}
\begin{proof}
By Proposition \ref{phicontinuous}, $\phi$ is continuous. To show that $\gamma$ is continuous,
we can assume that the bundles $\pi_{\xi_1, V}$
 and $\pi_{\xi_2, W}$ are trivial.
Based on (\ref{fp}),
we can express $\gamma$ as a continuous map from $W$
to the space of $r_2 \times r_1$ matrices, where $r_i$ is the rank of $\pi_{\xi_i, W}$, $i=1, 2$.
It follows from Lemma \ref{phifp} that $\phi$ is an $\tilde{\cf}$ map.
According to (\ref{fp}), $\gamma$ is an $\tilde{\cf}$ morphism.
\end{proof}

It is possible for the map $\phi$ in Proposition \ref{finiterank} to be constant on an open subset of $W$.


\section{Homomorphisms of mapping groups \label{mappinggroup}}

In this section, we apply the results from previous sections
to obtain a necessary and sufficient condition for a group homomorphism $f: \fvg \to \tilde{\cf}(W, G)$
to be a weighted composition operator. This condition can be further simplified
when $f$ is a group isomorphism $\fvg$ $\to$ $\fwg$.

For a Lie group $G$, we write $\autag$ (respectively  $\autg$) for the group of
algebraic group automorphisms (respectively  Lie group automorphisms) of $G$.
Let $S$ be a nonempty set and $\sg$ the space of all maps $S \to G$.

\begin{lemma} \label{pisa1.2}
  Let $f:$ $\fvg$ $\to$ $\sg$ be a group homomorphism.
  If $\dim_{\br} G \ge 2$, then the following statements (a) and (b) are equivalent.

  \begin{itemize}
    \item[(a)] For every $s \in S$, we have that
  $ (E_s \circ f)|_{\fczvg} \not\equiv \one$ and $ (E_s \circ f)|_{G}: G \to G$
  is surjective. Furthermore, the following condition holds:
  \begin{equation} \label{nonzero}
   f\left( \cf(V, G \setminus \{\one\}) \right) \subset \cm(S, G \setminus \{\one\}).
  \end{equation}

    \item[(b)]
  There exist maps $\phi: S \to V$ and $\gamma: S \times G \to G$
  such that $\gamma(s, \cdot) \in \autag$ for every $s \in S$ and
  \begin{equation} \label{gammaphi}
    f(x)(s)=\gamma(s, x \circ \phi(s)), \hspace{1mm} s \in S,
  \end{equation}
  for every $x \in \fvg$ with $x(V) \not=G$.
  \end{itemize}
    If $\dim_{\br} G \ge 2$, or if $\dim_{\br} G=1$ and $\cs_{E_s \circ f}$ consists of a single point
    for every $s \in S$, then statement (a) implies
  (\ref{gammaphi})
  for each $x \in \fczvg$. Additionally, if $G$ is a locally convex space $\sce$
  with $\dim_{\br} \sce \ge 2$, then statement (a) implies (\ref{gammaphi})
  for each $x \in \fve$.
\end{lemma}
\begin{proof}
  We only need to show that statement (a)  leads to the conclusions related to (\ref{gammaphi}). By (\ref{nonzero}),
  we have $ (E_s \circ f)|_{G} \in \autag$ for any $s \in S$.
  We define two maps $\gamma: S \times G \to G$ and $ h: \fvg \to \sg$ as follows:
  \begin{equation} \label{gammas}
  \gamma(s, \cdot)=(E_s \circ f)|_{G}, \hspace{1mm} s \in S, \hspace{2mm} \text{and}
  \end{equation}
  $$E_s \circ h(x)=\gamma(s, \cdot)^{-1} \circ E_s \circ f(x), \hspace{1mm} s \in S, \hspace{1mm} x \in \fvg,$$
  where $\gamma(s, \cdot)^{-1}$ is the inverse map of $\gamma(s, \cdot)$.
  Then $h$ is a group homomorphism such that $h(a)=a$ for any $a \in G$.
  Furthermore, we have $ (E_s \circ h)|_{\fczvg} \not\equiv \one$ for any $s \in S$,
  and $h$ still satisfies (\ref{nonzero}). Application of Proposition \ref{algebrard}
  along with Corollaries \ref{s1} and \ref{rdgh} to the map $h$ yields
  the conclusions related to (\ref{gammaphi}).
\end{proof}

In general, (\ref{gammaphi}) does not hold for all $x \in \fvg$ (Proposition \ref{su2}).

\begin{lemma} \label{fvghomo}
  Let $f: \fvg \to \tilde{\cf}(W, G)$ be a group homomorphism, and let
    $W_f \subset W$ be the open subset $\{w \in W: (E_w \circ f)|_{\fczvg} \not\equiv \one \}$.
    Assume the following conditions hold for each $w \in W_f$:
    The map $ (E_w \circ f)|_{G}: G \to G$
  is surjective. For every $x \in \cf(V, G \setminus \{\one\})$,
  we have $f(x)(w) \not=\one$.
  If $\dim_{\br} G=1$, then $\cs_{E_w \circ f}$ consists of a single point.
  Under these conditions,
  there exist
  a continuous map $\phi: W_f \to V$ and a map $\gamma: W_f \times G \to G$
  such that $\gamma(w, \cdot) \in \autag$ for every $w \in W_f$ and
    \begin{equation} \label{fx-2}
    f(x)(w)=\gamma(w, x \circ \phi(w)),  \hspace{1mm} w \in W_f,
  \end{equation}
 for
  every $x \in \fczvg$. When $\dim_{\br} G \ge 2$,
  the equation (\ref{fx-2})
  also holds for every $x \in \fvg$ with $x(V) \not=G$, and we have
  $ (E_w \circ f)|_{G} \equiv \one$ for every boundary point $w \in \partial W_f$.
\end{lemma}
\begin{proof}
Applying Lemma \ref{pisa1.2} to the group homomorphism $\fvg \ni x \mapsto f(x)|_{W_f} \in \cm(W_f, G)$
yields (\ref{fx-2}). The continuity of $\phi$ follows from
Proposition \ref{phicontinuous}.
 Let $w_0 \in \partial W_f$.
  Consider a pairwise distinct sequence $\{w_i\} \subset W_f$
  such that $\lim_{i \to \infty} w_i=w_0$.
By the definition of $W_f$ and (\ref{fx-2}), we have
  \begin{equation} \label{fxw0-2}
  \te \one=f(x)(w_0)=\lim_{i \to \infty}
  \gamma \left(w_i, x \circ \phi (w_i) \right), \hspace{2mm} x \in \fczvg.
  \end{equation}

  Assume for contradiction that the sequence $\{\phi(w_i)\}$ has no convergent subsequence.
  Suppose that $\cd$ is an open convex neighborhood of $\zero \in \fg$
such that $\exp_G|_{\cd}: \cd \to \exp_G(\cd)$ is a diffeomorphism,
$\{\phi(w_{i_k})\}$ is a pairwise distinct subsequence of $\{\phi(w_i)\}$ and $U$
is a second countable open neighborhood of the subsequence $\{\phi(w_{i_k})\}$.
Fix a metric on $U$. Choose
  pairwise disjoint
  open neighborhoods $V_k$ of $\phi(w_{i_k})$ for each $k \in \bn$
  such that
  $V_k \subset U$ and the diameter of $V_k$ is less than $k^{-1}$.
  Note that any sequence $\{v_k\}$ with $v_k \in V_k$ for each $k \in \bn$ does
  not have any convergent subsequence.
  Let $b_0 \in G \setminus \{\one\}$, and let $a_k$
  be the unique element of $G \setminus \{\one\}$ such that $\gamma(w_{i_k}, a_k)=b_0$.
  For each $k \in \bn$, we can find a curve $\Gamma_k \in C^{\infty}([0, 1], G)$ and
  a map $x_k \in C^{\infty}_c(V, G)$ such that
  $\supp x_k \subset V_k$, $x_k \circ \phi (w_{i_k})=\one$
  when $k$ is odd, $x_k \circ \phi (w_{i_k})=a_k$ when $k$ is even and $x_k(V) \subset \Gamma_k([0, 1])$.
  For any $v \in V$, there is at most one element $x_{k_v}$ in the sequence $\{x_k\}$
  such that $x_{k_v}(v) \not=\one$.
  Put
  \begin{equation} \label{x0}
  \te x_0=\prod_{k=1}^{\infty} x_{k} \in C^{\infty}(V, G).
  \end{equation}
  Note that $x_0(V) \not=G$ when $\dim_{\br} G \ge 2$.
  Set $x=x_0$ in (\ref{fx-2}). We obtain that
    $f(x_0)(w_0)=\lim_{k \to \infty} \gamma \left(w_{i_k}, x_k \circ \phi (w_{i_k}) \right)$.
    The limit on the right hand side does not exist.
  We have a contradiction.

  So the sequence $\{\phi(w_i)\}$ has a convergent subsequence
  $\phi(w_{i_j})$ $\to$ $v_0 \in V$.
  For any $a \in \exp_G(\cd)$, take
  $x_a \in \fczvg$ such that $x_a$ is constant $a$
  on a neighborhood of $v_0$.
  Set $x=x_a$ in (\ref{fxw0-2}) and choose $s=w_{i_j}$ in (\ref{gammas}). We have that
  $\one=\lim_{j \to \infty} \gamma \left(w_{i_j}, a \right)=\lim_{j \to \infty} f(a)(w_{i_j})=f(a)(w_0)$.
    Since the subset $\exp_G(\cd)$ generates the group $G$, we conclude that $ (E_{w_0} \circ f)|_{G} \equiv \one$.
\end{proof}

If $\partial W_f =\emptyset$, then
 $W_f$ must be the union of some of the components of $W$.
It is also possible for $\partial W_f \not=\emptyset$.
For example, let $n, m \in \bn$ with $m \ge 2$,
and let $\tau \in C^{\infty}_0(\br^n, \br)$ be a nonconstant function.
We define the linear map $f_1: C^{\infty}(\br^n, \br^m) \to C^{\infty}(\br^n, \br^m)$ by $f_1(x)=\tau x$.
This map $f_1$
satisfies all the conditions in Lemma \ref{fvghomo}
and we have $\partial W_{f_1} \not=\emptyset$.
We say that a bijection $\phi: W \to V$ is a bi-$\cf$ map, if both $\phi$ and $\phi^{-1}$
  are $\cf$ maps.

\begin{thm} \label{pisa1.2-iso}
  Let  $f:$ $\fvg$ $\to$ $\fwg$ be a group isomorphism.
  If $\dim_{\br} G \ge 2$, then the following statements (a) and (b) are equivalent.
  \begin{itemize}
    \item[(a)] The map
  $ (E_w \circ f)|_{G}: G \to G$
  is surjective for each $w \in W$ and
  \begin{equation} \label{nz=nz}
   f\left( \cf(V, G \setminus \{\one\}) \right) = \cf(W, G \setminus \{\one\}).
  \end{equation}

    \item[(b)] There exist a homeomorphism $\phi: W \to V$ and a map $\gamma: W \times G \to G$
  such that $\gamma(w, \cdot) \in \autag$ for each $w \in W$ and
  \begin{equation} \label{fxw}
    f(x)(w)=\gamma(w, x \circ \phi(w)), \hspace{2mm} w \in W,
  \end{equation}
  for any $x \in \fvg$ with $x(V) \not=G$.
  \end{itemize}
  If $G$ is a Banach Lie group, then statement (a) implies that $\gamma$ is continuous
  and $\phi$ is a  bi-$\cf$ map.
\end{thm}
\begin{proof}   It suffices to show that (a) implies (b).
By Lemma \ref{fvghomo}, there exist a nonempty open subset $W_f \subset W$,
  a continuous map $\phi: W_f \to V$ and a map $\gamma: W_f \times G \to G$
   such that $\gamma(w, \cdot) \in \autag$ for each $w \in W_f$, $\partial W_f=\emptyset$ and
  (\ref{fx-2}) holds
  for any $x \in \fvg$ with $x(V) \not=G$ and for
  any $x \in \fczvg$.
    We claim that
    \begin{equation} \label{v}
      V=\overline{\phi(W_f)}.
    \end{equation}
     Otherwise take a nonconstant map $x_0 \in \fczvg$ with  $\supp x_0 \subset V \setminus \overline{\phi(W_f)}$.
    It follows from the definition of $W_f$ and (\ref{fx-2}) that $f(x_0) \equiv \one$. Hence $f$ is not injective and we have a contradiction.
    Next we show that
     \begin{equation} \label{y}
    y(w)=\gamma(w, f^{-1}(y) \circ \phi(w)),  \hspace{1mm} w \in W_f,
  \end{equation}
  for any  $y \in \fwg$ with $y(W) \not=G$.
   Take $a \in G \setminus y(W)$. Then either $f^{-1}(a)=\one$ or
   $f^{-1}(a) \in \cf(V, G \setminus \{\one\})$ ((\ref{nz=nz})).
   Thus $f^{-1}(a)(V) \not=G$. Similarly, we have $f^{-1}(y a^{-1})(V) \not=G$.
   Choose $x=f^{-1}(a)$ and $x=f^{-1}(y a^{-1})$ in (\ref{fx-2}).
   Consequently, we obtain (\ref{y}).

   If $W \setminus W_f \not=\emptyset$,
  choose $y_0 \in \fwg \setminus \{\one\}$ such that $\supp y_0 \subset W \setminus W_f$ and $y_0(W) \not=G$.
  Note that $f^{-1}(y_0) \not\equiv \one$. By (\ref{v}) and (\ref{y}),
  $y_0|_{W_f} \not\equiv \one$. We have a contradiction. Thus $W_f=W$ and (\ref{fxw})
  holds
  for any $x \in \fvg$ with $x(V) \not=G$ and for
  any $x \in \fczvg$.

   We claim that $\phi$ is a bijection.
  Assume for contradiction that there exist distinct points $w_1, w_2 \in W$ such that $\phi(w_1)=\phi(w_2)$.
  Take $y_1 \in \fwg$ such that $y_1(W) \not=G$, $y_1(w_1)=\one$ and $y_1(w_2) \not=\one$.
  Set $y=y_1$ in (\ref{y}). This implies that either $y_1(w_1)=y_1(w_2)=\one$ or both $y_1(w_1)$ and $y_1(w_2)$
  are not equal to $\one$, leading to a contradiction. So $\phi$ is injective.
  Given $w_0 \in W$ and $y_i \in \fwg$ with $y_i(W) \not=G$, $i=2, 3$, it follows from (\ref{y})
   that
   \begin{equation} \label{w0}
   f^{-1}(y_2) \circ \phi(w_0)=f^{-1}(y_3) \circ \phi(w_0) \hspace{2mm}  \text{when} \hspace{2mm} y_2(w_0)=y_3(w_0).
   \end{equation}
  Suppose that $\phi$ is not surjective. Let $v^{\ast} \in V \setminus \phi(W)$.
  By (\ref{v}), there is a pairwise distinct sequence $\{w^{\ast}_j\}$ in $W$ such that
  $\phi(w^{\ast}_j) \to v^{\ast}$ as $j \to \infty$.
  Note that $\{w^{\ast}_j\}$ does not have any convergent subsequence (otherwise we would have $v^{\ast} \in \phi(W)$).
  Similar to the construction of the map in (\ref{x0}), we can find a map
  $y^{\ast} \in \fwg$ such that $y^{\ast}(W) \not=G$, $y^{\ast}(w^{\ast}_{2j-1})=\one$ and $y^{\ast}(w^{\ast}_{2j})=a^{\ast}$, where $j \in \bn$ and
  $a^{\ast}$ is a fixed element of  $G \setminus \{\one\}$. Note that $f^{-1}(a^{\ast}) \in \cf(V, G \setminus \{\one\})$.
  By (\ref{w0}),
  \begin{eqnarray*}
    & \te f^{-1}(y^{\ast})(v^{\ast})=\lim_{j \to \infty} f^{-1}(y^{\ast}) \circ \phi(w^{\ast}_{2j-1})=\one \hspace{1mm} \text{and} & \\
    & \te f^{-1}(y^{\ast})(v^{\ast})=\lim_{j \to \infty} f^{-1}(y^{\ast}) \circ \phi(w^{\ast}_{2j})=f^{-1}(a^{\ast})(v^{\ast}) \not=\one. &
  \end{eqnarray*}
  We have a contradiction.

  Let $a' \in G$, $v_0 \in V$ and
  $ a=f(a')(\phi^{-1}(v_0)) \in G.$
  By (\ref{w0}),
  $f^{-1}(a)(v_0)=a'$. Thus the map
  $ (E_{v_0} \circ f^{-1})|_{G}: G \to G$ is surjective, and $f^{-1}$ satisfies
  the same conditions as $f$. So there exist
  a continuous bijective map $\psi: V \to W$ and
  a map $\kappa: V \times G \to G$
  such that $\kappa(v, \cdot) \in \autag$ for any $v \in V$ and
    $f^{-1}(y)(v)=\kappa\left(v, y \circ \psi(v) \right)$ for any $v \in V$ and for any
    $y \in \fwg$ with
     $y(W) \not=G$.
  In view of (\ref{nz=nz}) and (\ref{fxw}), we have
$$x(v)=\kappa(v, f(x) \circ \psi(v))=\kappa(v, \cdot) \circ \gamma(\psi(v), \cdot) \circ x \circ \phi \circ \psi(v), \hspace{2mm} v \in V,$$
for any $x \in \cf(V, G \setminus \{\one\})$.
Set $x=a \in G \setminus \{\one\}$ in the equation above. We obtain that
 $\kappa(v, \cdot)=\gamma(\psi(v), \cdot)^{-1}$. Hence we have
 $x(v)= x \circ \phi \circ \psi(v)$, $v \in V$.
 With different choices of $x$,
  it follows that $\phi \circ \psi (v)=v$ for all $v \in V$.
Therefore $\psi=\phi^{-1}$
  and $\phi$ is a homeomorphism.

  If $G$ is a Banach Lie group, it follows from  Theorem \ref{frechetlg} that $\gamma$ is continuous
  and $\phi$ is a  bi-$\cf$ map.
  \end{proof}

  For a linear version of Theorem \ref{pisa1.2-iso}, see Theorem \ref{diffeo}.

\begin{cor} \label{bi-rd}
  Let $f: \fvg \to \fwg$ be a group
  isomorphism with
  \begin{equation} \label{fxw=xv}
  f(x)(W)=x(V), \hspace{1mm}  x \in \fvg.
  \end{equation}
  If $\dim_{\br} G \ge 2$,
  then there exists a bi-$\cf$ map $\phi: W \to V$ such that $ f(x)=x \circ \phi$
  for any $x \in \fvg$ with $x(V) \not= G$
  and for any $x \in \fczvg$.
\end{cor}
\begin{proof}
Note that $f(a)=a$, $a \in G$. The conclusion of the corollary follows
from Theorem \ref{pisa1.2-iso},
  equation (\ref{gammas}), Lemma \ref{pisa1.2} and Proposition \ref{phi}.
\end{proof}

It follows from Proposition \ref{su2} that the map $f$ in
Corollary  \ref{bi-rd} is not necessarily a composition operator
on $\fvg$.
If in Corollary \ref{bi-rd}, we set $V=W=S^1$, $G=\bc$ and
 $\cf=W^{1, p}$, where $p>1$, then $\phi$ is, in fact, a bi-Lipschitz map \cite[P. 710]{z17}.


\section{Linear maps between function spaces \label{lm}}

  Recall that $C^{\infty}(\br, \br)$ and $C^{\infty}(\br^2, \br)$
  are isomorphic as Fr\'echet spaces, despite the fact that $\br$ and $\br^2$ are not homeomorphic \cite[Theorem 5.3]{vo}.

  \begin{thm} \label{diffeo} Let  $n \in \bn$, and let
    $f:$ $\fvrn$ $\to$ $\fwrn$ be a linear bijection with
  \begin{equation} \label{rn}
   f( \cf(V, \br^n \setminus \{\zero \})) = \cf(W, \br^n \setminus \{\zero\}).
  \end{equation}
  Then there exist a bi-$\cf$ map $\phi: W \to V$
  and an $\cf$ map $\tilde{\gamma}: W \to GL(n, \br)$ such that
  \begin{equation} \label{8.1}
    f(x)(w)=\tilde{\gamma}(w) (x \circ \phi(w)), \hspace{2mm} w \in W, \hspace{2mm} x \in \fvrn.
    \end{equation}
  \end{thm}
\begin{proof} By (\ref{rn}), $(E_w \circ f)|_{\br^n}: \br^n \to \br^n$ is a linear bijection
for each $w \in W$. When $n \ge 2$, the conclusion of the theorem follows from
Theorem \ref{pisa1.2-iso} and Proposition \ref{finiterank}.

Now we consider the case when $n=1$.
We define $\tilde{\gamma}=f(1): W \to \br \setminus \{0\}$
and introduce a linear bijection $ h: \fvr \to \fwr$ given
by $h(x)=f(x)/\tilde{\gamma}$.
   It is clear that $h(a)=a$ for each $a \in \br$,
  and $h$ still satisfies (\ref{rn}). By Proposition \ref{algebrard},
  \begin{equation} \label{samerange}
   h(x)(W)=x(V), \hspace{2mm} x \in \fvr.
  \end{equation}
For any $v_0 \in V$, take $x_0 \in \fvr$
such that $ 0 \le x_0(v)<x_0(v_0)$ for each $v \in V \setminus \{v_0\}$. By (\ref{samerange}), $h(x_0)$ is non-negative and
there exists $w_0 \in W$ such that $h(x_0)(w_0)=x_0(v_0)$.
We claim that the function $h(x_0)$ attains its maximum value exclusively at $w_0$.
Otherwise we can find $w_1 \in W \setminus \{w_0\}$
with $h(x_0)(w_1)=x_0(v_0)$. Choose a cut-off function $\chi_0 \in C^{\infty}_c(W, [0, 1])$
such that $\chi_0$ is constant $1$ on a
neighborhood of $w_0$ and $w_1 \not\in \supp \chi_0$. Set $y_0=\chi_0 h(x_0)$ and $y_1=(1-\chi_0) h(x_0)$.
Note that $y_0(w_0)=x_0(v_0)=y_1(w_1)$. Both $h^{-1}(y_0)$ and $h^{-1}(y_1)$
are non-negative functions with maximal value $x_0(v_0)$. Suppose that
$h^{-1}(y_0)(v_1)=x_0(v_0)=h^{-1}(y_1)(v_2)$.
Since $x_0=h^{-1}(y_0)+h^{-1}(y_1)$,  it follows that $h^{-1}(y_1)(v_1)=0$. Thus $v_1 \not=v_2$,
and $x_0$
takes its maximal value at both $v_1$ and $v_2$. We have a contradiction.

Next we show that $\cs_{E_{w_0} \circ h}=\{v_0\}$.
For any open neighborhood $O$ of $v_0$, we denote by
$\chi_{O, v_0} \in C^{\infty}_c(O, [0, 1])$
a cut-off function that is constant $1$ on a neighborhood of $v_0$, where we consider
$C^{\infty}_c(O, [0, 1])$ as a subset of $C^{\infty}_c(V, [0, 1])$. If there exists $w_2 \in W \setminus \{w_0\}$
with $ h(\chi_{O, v_0} x_0)(w_2)=x_0(v_0)$, then $h(x_0)=h(\chi_{O, v_0} x_0)+h((1-\chi_{O, v_0}) x_0)$
also attains its maximal value $x_0(v_0)$ at $w_2 \not=w_0$. This leads to a contradiction.
Thus $ h(\chi_{O, v_0} x_0)(w_0)=x_0(v_0) >0$. Since $O$ is arbitrary, we have
$v_0 \in \cs_{E_{w_0} \circ h}$. For any $v_3 \in V \setminus \{v_0\}$,
take an open neighborhood $O_0$ of $v_0$ and an open neighborhood $O_3$ of $v_3$ such that
$\overline{O_0} \cap \overline{O_3}=\emptyset$. Choose a cut-off function $\chi_{O_0, v_0}$ described
as above. For any $x_1 \in \cf_c(O_3, \br)$, we can find a constant $c_1 \in \br \setminus \{0\}$
such that $h(c_1x_1)(w_0) \ge 0$ and $|c_1x_1(v)| <x_0(v_0)$ for any $v \in O_3$. Note that
the maximal value of $\chi_{O_0, v_0} x_0+c_1 x_1$
is $x_0(v_0)$ and $h(\chi_{O_0, v_0} x_0+c_1 x_1)(w_0)\ge x_0(v_0)$.
So we must have $h(x_1)(w_0)=0$. Hence $v_3 \not\in \cs_{E_{w_0} \circ h}$ and $\cs_{E_{w_0} \circ h}=\{v_0\}$.
Similarly $\cs_{E_{v_0} \circ h^{-1}}=\{w_0\}$.

From the above, we obtain a bijective map  $\phi: W \ni w_0 \mapsto v_0 \in V$.
It follows from (\ref{samerange}) and Lemma \ref{germ} that
\begin{equation} \label{h}
h(x)=x \circ \phi, \hspace{2mm} x \in \fcvr.
\end{equation}
By Proposition \ref{phi}, $\phi$ is a bi-$\cf$ map.

Let $x_2 \in \fvr$ be a function with an upper bound.
We claim
that if $x_2$ vanishes on a neighborhood $O_2$ of $v_0$, then $h(x_2)(w_0)=0$.
We may assume that $x_2(v)<x_0(v_0)$ for any $v \in V$ (otherwise replace $x_2$
by a constant multiple of it).
The maximal value of $h(\chi_{O_2, v_0} x_0+x_2)$ is $x_0(v_0)$ and $h(\chi_{O_2, v_0} x_0)(w_0)=x_0(v_0)$.
If $h(x_2)(w_0)\not=0$, then $h(x_2)(w_0)<0$.
Let $D_0$ be an open neighborhood of $w_0$ such that $\phi(D_0) \subset O_2$ and $h(x_2)$ is negative on $D_0$.
By (\ref{h}), we have $\supp h(\chi_{\phi(D_0), v_0} x_0) \subset D_0$. So
the maximal value of $h(\chi_{\phi(D_0), v_0} x_0+x_2)$ is strictly less than $x_0(v_0)$,
which contradicts (\ref{samerange}).

Let $x_3 \in \fvr$ be a function with an upper bound. Take an open neighborhood $O_4$ of $v_0$.
It follows from (\ref{h}) that
$$h(x_3)(w_0)=h(\chi_{O_4, v_0} x_3+(1-\chi_{O_4, v_0}) x_3)(w_0)=x_3 \circ \phi(w_0).$$
Fix a function $\tau \in C^{\infty}(\br, [0, 1])$ such that $\tau(t)=0$ for $t \ge 1$ and $\tau(t)=1$ for $t \le 0$.
For any $x \in \fvr$, we have $x=(\tau \circ x) x-(\tau \circ x-1) x$, which means $x$
is the difference of two $\cf$ functions that are bounded above. Hence $h(x)= x \circ \phi$,
and the conclusion of the theorem follows immediately.
 \end{proof}


\section{Homomorphisms of semigroups of functions \label{hsg}}

In this section, we apply results from Sections 6 and 7 to study
homomorphisms $f: \fvr \to \tfwr$
of multiplicative semigroups.

Let $\one \in \fvr$ (respectively $\zero \in \fvr$)
denote the constant function $1$ (respectively $0$), and let $w_0 \in W$. We say that $v_0 \in V$ is a support point
of the homomorphism $E_{w_0} \circ f$ of multiplicative semigroups
if, for every open neighborhood $O$ of $v_0$, there exists $x \in \fcvr$
such that $\supp x \subset O$ and $E_{w_0} \circ f(x) \not=0$.
We write $\cs_{E_{w_0} \circ f} \subset V$ for the subset of support points of $E_{w_0} \circ f$.
Note that $E_{w_0} \circ f(\one)$ (respectively $E_{w_0} \circ f(\zero)$)
can only take the values $1$ or $0$. If $E_{w_0} \circ f(\one)=0$ (respectively $E_{w_0} \circ f(\zero)=1$),
then $E_w \circ f(x)=0$ (respectively $E_w \circ f(x)=1$) for any $x \in \fvr$
and for any $w$ in the component of $W$ that contains $w_0$.

Let $\br^+$ be the group of positive real numbers, and let $x \in \fvrp$.
Note that $x=(\sqrt{x})^2$, where $\sqrt{x} \in \fvr$.
If $f(\one)=\one$, then $f^+=f|_{\fvrp}$ is a well-defined group homomorphism $\fvrp \to \tfwrp$.
In contrast to the case of group homomorphisms, the set
$\cs_{E_{w_0} \circ f}$ may be empty even if $E_{w_0} \circ f$ is surjective (cf. Proposition \ref{sp}).
For example, fix a point $v_1 \in V$ and define a homomorphism $f_{v_1}: \fvr \to \fvr$
of multiplicative semigroups as follows:
$f_{v_1}(x)=x$ if $x(v_1) \not=0$, and $f_{v_1}(x)=\zero$ if $x(v_1)=0$.
Then for any $v \in V \setminus \{v_1\}$, we have $\cs_{E_v \circ f_{v_1}}=\emptyset$
and $\cs_{E_v \circ f^+_{v_1}}=\{v\}$.

\begin{prop} \label{properties}
Let $f: \fvr \to \tfwr$ be a homomorphism of multiplicative semigroups, and let $w_0 \in W$.
    If $f(\one)=\one$, $f(\zero)=\zero$ and $v_0 \in \cs_{E_{w_0} \circ f}$, then
    the following statements hold.
    \begin{itemize}
      \item[(a)] $\cs_{E_{w_0} \circ f}=\{v_0\}$.

      \item[(b)] If $\xh_1$, $\xh_2 \in \fvr$ have the same germ at $v_0$ and $\xh_1(v_0) \not=0$, then
      $E_{w_0} \circ f(\xh_1)=E_{w_0} \circ f(\xh_2) \not=0$.
      If $\xh_3$ and $\zero$ have the same germ at $v_0$, then $E_{w_0} \circ f(\xh_3)=0$.

      \item[(c)] If in addition the set $E_{w_0} \circ f(\fvrp)$ contains a number other than $1$,
      then $\cs_{E_{w_0} \circ f^+}=\{v_0\}$.

      \item[(d)] If we further assume that there exists a map $\phi: W \to V$
      such that $\phi(w) \in \cs_{E_{w} \circ f}$
      for each $w \in W$, then $\phi$ is continuous.
    \end{itemize}
  \end{prop}
  \begin{proof} (a) Let $v_1 \in V \setminus \{v_0\}$. Take an open neighborhood
  $O_i$ of $v_i$ for $i=0, 1$, and choose $x_0 \in \fvr$ such that $O_0 \cap O_1=\emptyset$, $\supp x_0 \subset O_0$
  and $E_{w_0} \circ f(x_0) \not=0$. For any $x_1 \in \fvr$ with $\supp x_1 \subset O_1$,
  we have $x_0 x_1 = \zero$. Thus $E_{w_0} \circ f(x_1) =0$ and $v_1 \not\in \cs_{E_{w_0} \circ f}$.

(b) If $x_2 \in \fvr$ and $\one$ have the same germ at $v_0$,
take $x_3 \in \fvr$ such that $\supp x_3 \subset \{v \in V: x_2(v)=1 \}$
and $E_{w_0} \circ f(x_3) \not=0$. Then $f(x_2) f(x_3)=f(x_3)$,
which implies that $f(x_2)$ is $1$ on a neighborhood of $w_0$.
Now choose $x_4 \in \fvr$ such that $\xh_1 x_4$, $\xh_2 x_4$ and
$\one$ all have the same germ at $v_0$. Then we have
$E_{w_0} \circ f(\xh_1 x_4)=E_{w_0} \circ f(\xh_2 x_4)=1$.

Take an open neighborhood $O_{v_0}$ of $v_0$ and $x_5 \in \fcvr$ such that $\xh_3|_{O_{v_0}} \equiv 0$,
$\supp x_5 \subset O_{v_0}$ and $E_{w_0} \circ f(x_5) \not=0$. Then $\xh_3 x_5 = \zero$.
Hence $E_{w_0} \circ f(\xh_3)=0$.

(c) Take $x_6 \in \fvrp$ with $E_{w_0} \circ f(x_6)\not=1$.
For any open neighborhood $O$ of $v_0$, choose $x_7 \in \fvr$ such that $\supp x_7 \subset O$
and $x_7$ is constant $1$ on a neighborhood of $v_0$.
Then $e^{x_7 \ln x_6} \in \fvrp$ is constant $1$ on $V \setminus O$ and it shares the same germ at $v_0$ as $x_6$.
By (b), $E_{w_0} \circ f(e^{x_7 \ln x_6})=E_{w_0} \circ f(x_6)\not=1$. So $v_0 \in \cs_{E_{w_0} \circ f^+}$.
Let $v_1 \in V \setminus \{v_0\}$, and let $O_{v_1}$ be an open neighborhood of $v_1$
with $v_0 \not\in \overline{O_{v_1}}$.
For any $x_8 \in \fvrp$ with $\supp x_8 \subset O_{v_1}$
(note that the support of a Lie group valued function in $\fvrp$ differs from
the support of a continuous real valued function),
it follows from (b) that $E_{w_0} \circ f(x_8)=1$, which implies that $v_1 \not\in \cs_{E_{w_0} \circ f^+}$.

(d)
For any $w \in W$ and for any open neighborhood $O_{\phi(w)}$ of $\phi(w)$, choose
$x_9 \in \fvr$ and a neighborhood $D_w$ of $w$ such that $\supp x_9 \subset O_{\phi(w)}$ and
$E_{w'} \circ f(x_9) \not=0$ for each $w' \in D_w$. It follows from (b) that $\phi(D_w) \subset \supp x_9$.
Thus $\phi$ is continuous at $w$.
\end{proof}

\begin{prop} \label{isosp}
    Let $f: \fvr \to \fwr$ be an isomorphism of multiplicative semigroups.
    Then there exists a homeomorphism $\phi: W \to V$ such that
    $\cs_{E_w \circ f}=\{\phi(w)\}$ and $\cs_{E_{\phi(w)} \circ f^{-1}}=\{w\}$ for each $w \in W$.
\end{prop}
\begin{proof}
Take $\xh_1, \xh_2 \in \fvr$ such that $\xh_2 \not= \zero$ and
$\xh_1$
is constant $1$
(respectively $0$) on $\supp \xh_2$. Since $\xh_1 \xh_2 =\xh_2$ (respectively $\xh_1 \xh_2 = \zero$),
$f(\xh_1)$ is constant $1$ (respectively $0$) on the nonempty subset $\supp f(\xh_2) \subset W$.
In particular, if $\xh_3 \in \fvr$ is constant $1$ (respectively $0$) on some nonempty open subset of $V$,
then $f(\xh_3)$ is constant $1$ (respectively $0$) on some nonempty open subset of $W$.

 Let $w_0 \in W$. Choose a sequence $\{y_j\} \subset C^{\infty}_c(W, \br)$ such that
 $y_j(w_0)=1$, $y_{j}$ is constant $1$ on $\supp y_{j+1}$ for each $j \in \bn$, and
 for every open neighborhood $D$ of $w_0$, there exists $j_0 \in \bn$ such that
 $\supp y_j \subset D$ for all $j>j_0$.
 Assume for contradiction that $\cap_{j=1}^{\infty} \supp f^{-1}(y_{j})=\emptyset$.
 Note that  $f^{-1}(y_j)$
 is constant $1$ on $\supp f^{-1}(y_{j+1})$. Let $O_j \not=\emptyset$ be the
 interior of $\supp f^{-1}(y_{j}) \setminus \supp f^{-1}(y_{j+1})$.
 For each $j \in \bn$,
 take $x_j \in C^{\infty}_c(O_j, \br)$ such that $x_j$ is constant $1$
 on some open subset of $O_j$. Set $x^{\ast}=\sum_{k=1}^{\infty} x_{4k}$.
 Note that $x^{\ast}$ is a well-defined $C^{\infty}$ function on
 the complement of $\supp f^{-1}(y_{j})$ for each $j \in \bn$.
 Hence $x^{\ast} \in C^{\infty}(V, \br)$.
 Since $f(x^{\ast} x_{4k})=f(x_{4k}^2)$ and $f^{-1}(y_{4k-1}) x_{4k}=x_{4k}$,
 we have $f(x^{\ast})(w)=y_{4k-1}(w)=1$ if $f(x_{4k})(w)=1$.
 Similarly, $f(x^{\ast})(w)=0$ and $y_{4k-1}(w)=1$ if $f(x_{4k+2})(w)=1$.
 Thus for any open neighborhood $D$ of $w_0$, there exist points $w_{D, 1}, w_{D, 0} \in D$
 such that $f(x^{\ast})(w_{D, 1})=1$ and $f(x^{\ast})(w_{D, 0})=0$.
 This implies that $f(x^{\ast})$ is discontinuous at $w_0$, leading to a contradiction.
 So we conclude that $\cap_{j=1}^{\infty} \supp f^{-1}(y_{j}) \not=\emptyset$.

 Take $v_0 \in \cap_{j=1}^{\infty} \supp f^{-1}(y_{j})$. Since $f^{-1}(y_j)(v_0)=1$ for each $j \in \bn$,
 we have $w_0 \in \cs_{E_{v_0} \circ f^{-1}}$. Similarly, $v_0 \in \cs_{E_{w_1} \circ f}$ for some $w_1 \in W$.
 If $w_1 \not=w_0$,
 it is clear that $y_j(w_1)=0$  for sufficiently large $j$. Note that $f^{-1}(y_j)$
 and $\one$ have the same germ at $v_0$.
 It follows from Proposition \ref{properties}(b) that $v_0 \not\in \cs_{E_{w_1} \circ f}$. We have a contradiction.
 So $w_1=w_0$. By Proposition \ref{properties}(a), $\cs_{E_{w_0} \circ f}=\{v_0\}$
 and $\cs_{E_{v_0} \circ f^{-1}}=\{w_0\}$.
 Let $\phi$ be the map $W \ni w_0 \mapsto v_0 \in V$.
 According to Proposition \ref{properties}(d), $\phi$ is a homeomorphism.
\end{proof}

\begin{thm} \label{smgp}
  Let $f: \fvr \to \tfwr$ be a homomorphism of multiplicative semigroups. If
  $f(\cf(V, \br^+ \setminus \{1\})) \subset \tilde{\cf}(W, \br^+ \setminus \{1\})$,
  $E_w \circ f(\br^+)=\br^+$ for
 each $w \in W$ and there exists a map $\phi: W \to V$ such that $\phi$ is not constant on any open subset of $W$
 and $\phi(w) \in \cs_{E_w \circ f}$ for each $w \in W$,
 then the following statements hold.
 \begin{itemize}
   \item[(a)] There exists an $\tilde{\cf}$  function $\alpha: W \to \br^+$
    such that for each $w \in W$ and for each $x \in \fvr$ with $x \circ \phi(w) \not=0$, we have
    \begin{equation} \label{fplus}
    \te f(x)(w)=f(\text{sgn}(x \circ \phi(w))) |x \circ \phi(w)|^{\alpha(w)},
    \end{equation}
    where we consider $\text{sgn}(x \circ \phi(w))$ as a constant function on $V$.
    Furthermore, $\phi$ is an $\tilde{\cf}$ map.

   \item[(b)] If in (a), $\phi$ is an open map,
    then (\ref{fplus}) holds for each $w \in W$ and for each
     $x \in \fvr$.

   \item[(c)]  If in (a), $\fvr \subset C^1(V, \br)$ and
   $f(x) \subset C^1(W, \br)$ for each $x \in \fvr$, then $\alpha(w) \ge 1$ for each $w \in W$.
 \end{itemize}
\end{thm}
\begin{proof}
  (a) It is clear that $f(\one)=\one$ and $f(\zero)=\zero$.
  According to Proposition \ref{properties}(c), $\cs_{E_w \circ f^+}=\{\phi(w)\}$
  for each $w \in W$.
  Any nonconstant Lie group homomorphism $\br^+ \to \br^+$
  is of the form $t \mapsto t^{\alpha}$, where $\alpha \in \br \setminus \{0\}$
  is a constant.
  Application of  Lemma \ref{fvghomo} and Theorem \ref{frechetlg} to $f^+$
   leads to the existence of a map $\alpha: W \to \br \setminus \{0\}$, ensuring that
  (\ref{fplus}) holds for all $x \in \fczvrp=\cf_c(V, \br^+)$.
  By  Theorems \ref{frechetlg} and \ref{wco2},
  $\phi$ and $\alpha$ are $\tilde{\cf}$ maps.
  Let $x^{\ast} \in \fvr$ and $w^{\ast} \in W$.
  If $x^{\ast} \circ \phi(w^{\ast})>0$, take a positive number $\ve$
  and a function $\tau \in C^{\infty}_c(V, [0, 1])$ such that $\ve < x^{\ast} \circ \phi(w^{\ast})$,
  $\tau$ is constant $1$ on a neighborhood of $\phi(w^{\ast})$ and
  $\supp \tau \subset \{v \in V: x^{\ast}(v)>\ve \}$. Then the function
  $e^{\tau \ln(\tau(x^{\ast}-\ve)+\ve)} \in \cf_c(V, \br^+)$
   shares the same germ at $\phi(w^{\ast})$ as $x^{\ast}$.
  It follows from Proposition \ref{properties}(b) that
  (\ref{fplus})  is valid when $x=x^{\ast}$ and $w=w^{\ast}$.
  If $x^{\ast} \circ \phi(w^{\ast})<0$, then $f(x^{\ast})=f(-1)f(-x^{\ast})$,
  and (\ref{fplus}) remains valid when $x=x^{\ast}$ and $w=w^{\ast}$.

 Next we show that $\alpha$ is a positive function.
  Assume for contradiction that there exists $w_0 \in W$ such that
  $\alpha$ is negative on the component of $W$ containing $w_0$.
  Take a sequence $\{w_k\} \subset W \setminus \{w_0\}$ such that
  $w_k \to w_0$ and the points $\phi(w_0), \phi(w_1), \cdots$ are pairwise distinct. Choose $x_0 \in C^{\infty}(V, \br)$
  such that $x_0 \circ \phi(w_k) >0$ for each $k \in \bn$ and $x_0 \circ \phi(w_0)=0$ (Lemma \ref{passsequence}).
  Set $x=x_0$, $w=w_k$ in (\ref{fplus})
  and take the limits on both sides as $k \to \infty$.  This leads to a contradiction.

  (b)
  We claim that if $\phi$ is an open map and $x \circ \phi(w_0)=0$,
  then $f(x)(w_0)=0$. If $x$ is constant $0$
  on a neighborhood of $\phi(w_0)$, this follows from Proposition \ref{properties}(b).
  If $x$ is not constant $0$ on
  any neighborhood of $\phi(w_0)$, choose precompact open neighborhoods $D_k$, $k \in \bn$, of $w_0$ such that
  $\overline{D_{k+1}} \subset D_k$ and $\cap_{k=1}^{\infty} D_k=\{w_0\}$.
  Note that each $\phi(D_k)$ is an open neighborhood of $\phi(w_0)$.
  Take $w_k \in D_k$ with $x \circ \phi(w_k) \not=0$, $k \in \bn$. Then $w_k \to w_0$.
  Set $w=w_k$ in (\ref{fplus}) and take the limits as $k \to \infty$. This yields $f(x)(w_0)=0$.

  (c) In this case, $\phi$ is a $C^1$ map.
  We only
  need to show that if the differential of $\phi$ at $w_0^{\ast} \in W$ is not zero,
  then $\alpha(w_0^{\ast}) \ge 1$. Take a $C^{\infty}$ curve $w(t): (-1, 1) \to W$
  and $x_0^{\ast} \in C^{\infty}(V, \br)$
  such that $w(0)=w_0^{\ast}$, $x_0^{\ast} \circ \phi(w_0^{\ast})=0$
  and the differential of $x_0^{\ast} \circ \phi \circ w(t)$ at $t=0$ is not zero.
  Then there exists a neighborhood $U \subset (-1, 1)$ of $0$ such that
  $x_0^{\ast} \circ \phi \circ w(t) \not=0$ for every $t \in U \setminus \{0\}$.
  Set $x=x_0^{\ast}$, $w=w(t)$ in  (\ref{fplus}) and take the limits as $t \to 0$. This gives
  that $f(x_0^{\ast})(w_0^{\ast})=0$.
  It is straightforward to verify that if $0<\alpha(w_0^{\ast})<1$, then the right hand side of (\ref{fplus})
  is not differentiable at $t=0$. This leads to a contradiction.
\end{proof}

If the map $\phi$ in Theorem \ref{smgp} is constant, then $\alpha$
might be negative.  Fix a point $v_0 \in V$. We define $f(x) = 1/x(v_0)$ if $x(v_0) \not=0$,
and $f(x) = \zero$ if $x(v_0) =0$. In this case, $\phi(w) = v_0$ and $\alpha(w) =-1$ for each $w \in W$.

\begin{prop} \label{smallerthan1}
  Let $V_0$ be an open subset of $\br^p$, where $p \in \bn$,
  $f: C^k(V_0, \br) \to C(W, \br)$, where $k=0, 1, \cdots, \infty$, a homomorphism of multiplicative semigroups,
   $w_0 \in W$,
  $v_0 \in \cs_{E_{w_0} \circ f}$ and $O$ a
  connected open neighborhood of $v_0$.
  If there  exist
  $x_0 \in C^k(V_0, \br)$ and $\tau \in C^{\infty}_c(O, [0, 1])$ such that
  $\tau$ is constant $1$
  on a neighborhood of $v_0$, $f(x_0)(w_0)>1$ and $0 \le x_0(v)<c$ for all $v \in \supp \tau$, where $c<1$ is a constant,
  then for each $j \in \bn$, there exists $n_j \in \bn$  such that
  $f(\tau x_0^{n_j})(w_0)>j$ and the $C^k$ norm
  $||\tau x_0^{n_j}||_{C^k}<2^{-j}$ when $k <\infty$ (respectively the $C^j$ norm
   $||\tau x_0^{n_j}||_{C^j}<2^{-j}$ when $k =\infty$).
\end{prop}
\begin{proof}
  The condition $f(x_0)(w_0)>1$ implies that $f(\one)=\one$ and $f(\zero)=\zero$
  on the component of $W$ containing $w_0$.
  By Proposition \ref{properties}(b), $ f(\tau x_0^n)(w_0)=(f(x_0)(w_0))^n \to \infty$ when $n \to \infty$.
  It suffices to show that for any fixed multi-index $\alpha_0=(\alpha^1_0, \cdots, \alpha^p_0)$,
  where $a^1_0$, $\cdots$, $\alpha^p_0 \in \bn \cup \{0\}$, such that
  $|\alpha_0|=\alpha_0^1+\cdots+\alpha^p_0 \le k$ when $k$ is finite and $|\alpha_0| \le j$ when $k=\infty$,
  the sequence of functions $\{\partial^{\alpha_0} (\tau x_0^n)\}$
  converges uniformly to $\zero$ on $O$ as $n \to \infty$.
 We say that $\alpha=(\alpha^1, \cdots, \alpha^p) \le \alpha_0$ if $\alpha^i \le \alpha_0^i$ for all $i=1, \cdots, p$.
 Let $M>c$ be a constant such that
 $$\te \sup_{\alpha \le \alpha_0, v \in \supp \tau} \left|\partial^{\alpha} x_0(v) \right|, \hspace{1mm}
 \sup_{\alpha \le \alpha_0, v \in \supp \tau} \left|\partial^{\alpha} \tau(v) \right| \le M.$$
 The function $\partial^{\alpha_0} (\tau x_0^n)$
 is the sum of $(n+1)^{|\alpha_0|}$ terms. Each term
 is a product of $n+1$ factors, where each factor
 is either a derivative of $x_0$ or  a derivative of $\tau$. When $n>|\alpha_0|$, each term is divisible by $x_0^{n-|\alpha_0|}$.
 Therefore
$ |\partial^{\alpha_0} (\tau x_0^n)| \le (n+1)^{|\alpha_0|} c^{n-|\alpha_0|} M^{|\alpha_0|+1} \to 0$
as $n \to \infty$.
\end{proof}

\begin{thm} \label{semigp}
Let $f: \ckovr \to \cktwr$, where $k_1, k_2 =0, 1, \cdots, \infty$
and $k_1 \ge k_2$, be a homomorphism of multiplicative semigroups.
If $E_w \circ f|_{\ckovrp}$ is not constant for any $w \in W$ and there exists
an open map $\phi: W \to V$
   such that $\phi(w) \in \cs_{E_w \circ f}$ for any $w \in W$,
   then $\phi$ is a $C^{k_2}$ map and there exists a $C^{k_2}$ map $\alpha: W \to \br^+$ such that
\begin{equation*}
      f(x)(w)=f(\mathrm{sgn}(x \circ \phi(w))) |x \circ \phi(w)|^{\alpha(w)}, \hspace{2mm} w \in W, \hspace{2mm} x \in \ckovr.
    \end{equation*}
If $k_2 \ge 1$, then $\alpha(w) \ge 1$ for any $w \in W$.
\end{thm}
\begin{proof}
It is clear that $f(\one)=\one$ and $f(\zero)=\zero$.
By Proposition \ref{properties}(d),
$\phi$ is continuous.
We claim that if $0<x_0 \circ \phi(w_0)<1$, where $w_0 \in W$ and $x_0 \in \ckovr$, then $0<f(x_0)(w_0) \le 1$.
Assume, for the sake of contradiction,  that $f(x_0)(w_0)>1$. Take
a precompact open neighborhood $D$ of $w_0$ and a local chart $(U, \Phi)$ of $V$
 such that  $\phi(\overline{D}) \subset U$ and for any $w \in D$, we have
$0<x_0 \circ \phi(w) <1$ and $f(x_0)(w)>1$.
Choose a sequence $\{w_j\} \subset D \setminus \{w_0\}$
such that $w_j \to w_0$ and the points $\phi(w_0), \phi(w_1), \phi(w_2), \cdots$ are
pairwise distinct. For each $j \in \bn$, consider a precompact open neighborhood $O_j$ of $\phi(w_j)$
and a function $\tau_j \in C^{\infty}_c(O_j, [0, 1])$ such that
$\overline{O_j} \subset \phi(D)$, $\phi(w_0) \not\in \overline{O_j}$,
$\overline{O_{j_1}} \cap \overline{O_{j_2}} =\emptyset$ for $j_1 \not=j_2$ and
$\tau_j$ is constant $1$
on a neighborhood of $\phi(w_j)$.
By Proposition \ref{smallerthan1}, for each $j \in \bn$, there exists
$n_j \in \bn$ such that
$f(\tau_j x_0^{n_j})(w_j)>j$ and  $||(\tau_j x_0^{n_j}) \circ \Phi^{-1}||_{C^{k_1}}<2^{-j}$
when $k_1<\infty$ (respectively $||(\tau_j x_0^{n_j}) \circ \Phi^{-1}||_{C^{j}}<2^{-j}$ when $k_1=\infty$).
Note that $x^{\ast}=\sum_{j=1}^{\infty} \tau_j x_0^{n_j}$
is a well-defined $C^{k_1}$ function on $V$ with $\supp x^{\ast} \subset \phi(\overline{D})$.
By Proposition \ref{properties}(b), $f(x^{\ast})(w_j)>j$. Thus the
continuous function $f(x^{\ast})$ is unbounded on the compact subset $\overline{D}$,
 which leads to  a contradiction.

Next we show that the map $E_{w_0} \circ f|_{\br^+}$ is not constant.
Otherwise $E_{w_0} \circ f|_{\br^+} \equiv 1$.
For any $x \in \ckovrp$, we can find $t \in \br^+$ such that $tx \circ \phi(w_0)<1$. Thus $f(x)(w_0) \le 1$.
Similarly $f(x^{-1})(w_0) \le 1$. So $f(x)(w_0) \ge 1$, and $E_{w_0} \circ f|_{\ckovrp}$ is constant, which
leads to a contradiction.
If $x_1 \circ \phi(w_0) <x_2 \circ \phi(w_0)$, where $x_1, x_2 \in C^{k_1}(V, \br^+)$,
then $0<(x_1 x_2^{-1}) \circ \phi(w_0)<1$.
Hence $f(x_1)(w_0) \le f(x_2)(w_0)$. In particular, $E_{w_0} \circ f|_{\br^+}$
is an increasing nonconstant group homomorphism. So $E_{w_0} \circ f(t)=t^{\alpha(w_0)}$,
where $t \in \br^+$ and $\alpha(w_0)>0$. We claim that the condition $0<x_0 \circ \phi(w_0)<1$
implies that $0<f(x_0)(w_0) < 1$. Suppose $f(x_0)(w_0) = 1$.
Then there exists $t>1$ such that $0<t x_0 \circ \phi(w_0)<1$
and $f(t x_0)(w_0)=t^{\alpha(w_0)}>1$. This again leads to a
contradiction. Therefore $f(C^{k_1}(V, \br^+ \setminus \{1\})) \subset C^{k_2}(W, \br^+ \setminus \{1\})$.
The conclusion of the theorem follows from  Theorem \ref{smgp}.
\end{proof}

\begin{cor} \label{semigpiso}
  Let $f: \ckvr \to \ckwr$, where
  $k=0, 1, \cdots, \infty$,
  be an isomorphism of multiplicative semigroups. If $k \ge 1$, then there exists a $C^k$ diffeomorphism $\phi: W \to V$
  such that
  $$f(x)=x \circ \phi, \hspace{2mm} x \in \ckvr.$$
  If $k=0$, then there exist a homeomorphism $\phi: W \to V$
  and a continuous positive function $\alpha$ on $W$ such that
  \begin{equation*}
      f(x)(w)=\mathrm{sgn}(x \circ \phi(w)) |x \circ \phi(w)|^{\alpha(w)}, \hspace{2mm} w \in W, \hspace{2mm} x \in C(V, \br).
    \end{equation*}
\end{cor}
\begin{proof}
  The conclusion of the corollary follows from Proposition \ref{isosp} and Theorem \ref{semigp}.
\end{proof}


\section{Homomorphisms of groups of sections of bundles \label{hlgb}}

Let $\pi_{\scg_1, V}: \scg_1 \to V$ and $\pi_{\scg_2, W}: \scg_2 \to W$ be $C^{\infty}$ Lie group
bundles with the same typical fiber $G$. We denote the evaluation map $\fvscgo \ni x \mapsto x(v) \in (\scg_1)_v$ at $v \in V$
by $E_v$.

\begin{thm} \label{lgb}
 Let $f: \fczvscgo \to \ftscgt$ be a group homomorphism.
 Then the following statements (a) and (b) are equivalent.
 \begin{itemize}
   \item[(a)]  For any $w \in W$, we have $\cs_{E_w \circ f} \not=\emptyset$,
   and there exist open subsets $O_w$, $O'_w \subset V$,
   a local trivialization $\Psi_w: \pi_{\scg_1, V}^{-1}(O'_w) \to O'_w \times G$ of $\pi_{\scg_1, V}$
   and a map $\varrho_w: G \to \cf_c^0(\pi_{\scg_1, V}|_{O'_w}) \subset \fczvscgo$
   such that $O_w$ is precompact, $\overline{O_w} \subset O'_w$,
   $\cs_{E_w \circ f} \subset O_w$,
   $E_w \circ f(x) \not=\one$ for any $x \in \cf_c^0(\pi_{\scg_1, V}|_{O'_w})$
   with  $x(v) \not=\one$, where $v \in O_w$,
   $\Psi_w \circ \varrho_w(a)|_{O_w}$ is the constant section $a \in G$  and
   $E_w \circ f \circ \varrho_w: G \to (\scg_2)_w$
   is a surjective group homomorphism.
   If $\dim_{\br} G=1$, we further require that $\cs_{E_w \circ f}$ consists of a single point
   for any $w \in W$.

   \item[(b)] There exist a continuous map $\phi: W \to V$ and a (not necessarily continuous) bundle morphism
   $\gamma: \phi^{\ast} \scg_1 \to \scg_2$ over the identity map on $W$
   such that the restriction of $\gamma$ to each fiber
   is a group isomorphism $(\scg_1)_{\phi(w)} \to (\scg_2)_w$ and
    \begin{equation} \label{fcxo}
      f(x)=\gamma \circ (\phi^{\ast} x), \hspace{2mm} x \in \fczvscgo.
    \end{equation}
 \end{itemize}
 If $\pi_{\scg_1, V}$ and $\pi_{\scg_2, W}$ are Banach  Lie group
 bundles and
  $\phi$ is not constant on any open subset of $W$,
 then $\phi$ is an $\tilde{\cf}$ map and
  $\gamma$ is continuous.
\end{thm}
\begin{proof}
First we show that (b) implies (a). Note that $\cs_{E_w \circ f}=\{\phi(w)\}$.
Take precompact open neighborhoods $O_w$, $O'_w$ of $\phi(w)$
such that $\overline{O_w} \subset O'_w$, there exists
a local trivialization $\Psi_w: \pi^{-1}_{\scg_1, V}(O'_w) \to O'_w \times G$
of $\pi_{\scg_1, V}$ and we can find a function $\chi \in C^{\infty}_c(O'_w, [0, 1])$
with $\chi|_{O_w} \equiv 1$.
Recall that any open neighborhood of $\one \in G$ generates the group $G$.
Hence for any $a \in G$, we can find finitely many points $\ah_1, \cdots, \ah_p$ in the Lie algebra $\fg$ of $G$
such that $a=\exp_G(\ah_1) \cdots \exp_G(\ah_p)$. Set $\xh_j=\chi \ah_j \in C^{\infty}_c(O'_w, \fg)$, $j=1, \cdots, p$.
Define $x_a \in \cf_c^0(\pi_{\scg_1, V}|_{O'_w})$ by
$x_a(v)= \Psi_w^{-1}(v, \exp_{G} \xh_1(v) \cdots \exp_{G} \xh_p(v))$, where $v \in O'_w$,
and define a map $\varrho_w: G \to \cf_c^0(\pi_{\scg_1, V}|_{O'_w})$ by
$\varrho_w(a)=x_a$, where $a \in G$.
It is straightforward to verify that (a) holds.

Next we show that (a) implies (b).
For each $y \in \cf^0_c(O'_w, G)$, we write $\Psi_w^{-1} y \in \cf_c^0(\pi_{\scg_1, V}|_{O'_w})$ for the section
$O'_w \ni v \mapsto \Psi_w^{-1}((v, y(v))) \in \pi^{-1}_{\scg_1, V}(O'_w)$.
Note that $E_w \circ f \circ \varrho_w: G \to (\scg_2)_w$ is a group isomorphism.
Define a group homomorphism by
\begin{equation} \label{gw}
g_w: \cf^0_c(O'_w, G) \ni y \mapsto (E_w \circ f \circ \varrho_w)^{-1} \circ E_w  \circ f(\Psi_w^{-1} y) \in G.
\end{equation}
Then $g_w(y) \not=\one$ for any $y \in \cf^0_c(O'_w, G)$
   such that  $y(v) \not=\one$, where $v \in O_w$.
Let $P_2: O'_w \times G \to G$ be the projection onto the second component.
Set $y_a=P_2 \circ \Psi_w \circ \varrho_w(a) \in \cf^0_c(O'_w, G)$, where $a \in G$.
Then $y_a|_{O_w} \equiv a$ and $g_w(y_a)=a$ for each $a \in G$.
For any $y \in \cf^0_c(O'_w, G)$, we have $ g_w(y y_{g_w(y)^{-1}})=\one$.
Thus there exists $v_y \in O_w$ such that $y y_{g_w(y)^{-1}} (v_y)=\one$, which implies that $y(v_y)=g_w(y)$.
Hence the map $g_w|_{\cf^0_c(O_w, G)}$ is range decreasing.

Since $\cs_{E_w \circ f} \subset O_w$, we have $g_w|_{\cf^0_c(O_w, G)} \not\equiv \one$.
According to Theorem \ref{nbhd1}, there exists $\vb(w) \in O_w$ such that
\begin{equation} \label{gwy}
g_w(y)=y(\vb(w)), \hspace{2mm} y \in \cf^0_c(O_w, G).
\end{equation}
Note that $\cs_{E_w \circ f}=\{\vb(w)\}$.
Define the map $\phi: W \to V$ by $\phi(w)=\vb(w)$. For an element $a' \in (\scg_1)_{\phi(w)}$,
if $\Psi_w(a')=(\phi(w), a)$ for some $a \in G$,
we denote $a$ as $\Psi_{w, \phi(w)}(a')$.
Set $y=P_2 \circ \Psi_w \circ x$ in (\ref{gwy}), where $x \in \cf_c^0(\pi_{\scg_1, V}|_{O_w})$.
In view of (\ref{gw}), we obtain that
\begin{equation*}
E_w \circ  f(x)=(E_w \circ f \circ \varrho_w \circ \Psi_{w, \phi(w)})(x \circ \phi(w)),
\hspace{2mm} x \in \cf_c^0(\pi_{\scg_1, V}|_{O_w}).
\end{equation*}
By Proposition \ref{sp}, the equation above also holds for any $x \in \fczvscgo$.
Define the map $\gamma: \phi^{\ast} \scg_1 \to \scg_2$ by $\gamma|_{(\scg_1)_{\phi(w)}}
=E_w \circ f \circ \varrho_w \circ \Psi_{w, \phi(w)}: (\scg_1)_{\phi(w)}
\to (\scg_2)_w$.
Then we have (\ref{fcxo}).

By Proposition \ref{phicontinuous},
$\phi$ is continuous.
 If $\pi_{\scg_1, V}$ and $\pi_{\scg_2, W}$ are Banach  Lie group
 bundles and
  $\phi$ is not constant on any open subset of $W$,
  it
  follows from Theorem \ref{frechetlg} that $\phi$ is an $\tilde{\cf}$ map and
  $\gamma$ is continuous.
\end{proof}

\begin{thm} \label{finiterankbundle}
  Let $\pi_{\eta_1, V}: \eta_1 \to V$ and $\pi_{\eta_2, W}: \eta_2 \to W$ be $C^{\infty}$ real vector bundles of the same rank $r \in \bn$,
  and let $f: \foceo \to \ftet$ be a linear map.
 Then the following statements (a) and (b) are equivalent.
 \begin{itemize}
   \item[(a)] For any $w \in W$, $\cs_{E_{w} \circ f}$ is a nonempty and compact
   subset of $V$. Moreover,
   the bundle $\pi_{\eta_1, V}$ is trivial on an open neighborhood of $\cs_{E_{w} \circ f}$.
   Additionally,
   for any $x \in \foceo$ which does not vanish on $\cs_{E_{w} \circ f}$, we have $f(x)(w) \not=\zero$.
   If $r=1$, we further require that $\cs_{E_w \circ f}$ consists of a single point
   for every $w \in W$.

   \item[(b)] There exist
  an $\tilde{\cf}$ map $\phi: W \to V$ and an $\tilde{\cf}$ isomorphism of vector bundles
  $\gamma: \phi^{\ast} \eta_1 \to \eta_2$
  covering the identity map on $W$
  such that
  \begin{equation} \label{gwcovb}
  f(x)=\gamma \circ (\phi^{\ast} x), \hspace{2mm} x \in \foceo.
  \end{equation}
 \end{itemize}
\end{thm}
\begin{proof} One direction being trivial, we shall only verify that
(a) implies (b).
Take precompact open neighborhoods $O_w$, $O'_w$ of $\cs_{E_w \circ f}$ such that $\overline{O_w} \subset O'_w$,
there exists a local trivialization $\Psi_w: \pi^{-1}_{\eta_1, V}(O'_w) \to O'_w \times \br^r$ of $\pi_{\eta_1, V}$
and for any $a \in \br^r$, we can find $x_a \in \foceo$ such that $\supp x_a \subset O'_w$ and
$\Psi_w \circ x_a|_{O_w}$ is the constant section $a$.
Define a map $\varrho_w: \br^r \to \foceo$ by $\varrho_w(a)=x_a$. According to Proposition \ref{sp},
$E_w \circ f \circ \varrho_w: \br^r \to (\eta_2)_w$ is a well-defined
linear map (even though $\varrho_w$ is not necessarily linear). Note that $E_w \circ f \circ \varrho_w$
is actually an isomorphism of vector spaces. The conclusion of the theorem
follows from Theorem \ref{lgb}
 and Proposition \ref{finiterank}.
\end{proof}


\section{Homomorphisms of Lie algebras of vector fields \label{vf}}

Let $\citv$ (respectively $\xv$) be the Lie algebra of all $C^{\infty}$ vector fields
   (respectively of all $C^{\infty}$ vector fields with compact support) on $V$.
   The following theorem can be considered as a homomorphism version of the Shanks-Pursell theorem \cite{sp}.

\begin{thm} \label{shanks}
   Suppose that $\dim_{\br} V=\dim_{\br} W$ and $f: \xv \to \citw$ is
   a Lie algebra homomorphism. Then $f$  satisfies the conditions in
   Theorem \ref{finiterankbundle}(a) if and only if
 there exists a $C^{\infty}$ local diffeomorphism $\phi: W \to V$
   such that
   \begin{equation} \label{dphi}
   f(x)(w)=(d_{w}\phi)^{-1}(x \circ \phi(w)), \hspace{2mm} w \in W, \hspace{2mm} x \in \xv,
   \end{equation}
   where $(d_{w}\phi)^{-1}$
   is the inverse of the differential of $\phi$ at $w$.
\end{thm}
\begin{proof}
  It is straightforward to verify  the ``if" direction.   In the ``only if"
direction,
  it follows from Theorem \ref{finiterankbundle} that there exist a $C^{\infty}$ map
  $\phi: W \to V$ and a $C^{\infty}$ isomorphism
  of vector bundles $\gamma: \phi^{\ast} TV \to TW$
  covering the identity map on  $W$ such that
  \begin{equation} \label{wcomposition}
  f(x)=\gamma \circ (\phi^{\ast} x), \hspace{2mm} x \in \xv.
  \end{equation}

  For any $w_0 \in W$, choose a precompact open neighborhood $D_{w_0}$ of $w_0$ and open neighborhoods
  $O_{\phi(w_0)}$, $O'_{\phi(w_0)}$ of $\phi(w_0)$ such that
  $O_{\phi(w_0)}$ is precompact,
  $\overline{O_{\phi(w_0)}} \subset O'_{\phi(w_0)}$,
  $\phi(\overline{D_{w_0}}) \subset O_{\phi(w_0)}$,
$O'_{\phi(w_0)}$ is contained in a local chart of $V$
and there exists a function $\chi_{\phi(w_0)} \in C^{\infty}_0(O'_{\phi(w_0)}, [0, 1])$
with $\chi_{\phi(w_0)}|_{O_{\phi(w_0)}} \equiv 1$.
Let $(v^1, \cdots, v^p)$, where $p=\dim_{\br} V$, be local coordinates on
$O'_{\phi(w_0)}$.
Set $X_j=\chi_{\phi(w_0)} \partial/\partial v^j \in \xv$, $j=1, \cdots, p$.
It follows from (\ref{wcomposition}) that
$f(X_1)(w), \cdots, f(X_p)(w)$ form a basis for $T_w W$ when $w \in D_{w_0}$.
For any $x' \in \xv$,
there exists $x \in \xv$ such that $\supp x \subset O_{\phi(w_0)}$ and $f(x)|_{D_{w_0}}=f(x')|_{D_{w_0}}$.
We can write $x=\sum_{j=1}^p \alpha_j X_j$, where $\alpha_j \in C^{\infty}_c(O_{\phi(w_0)})$.
By (\ref{wcomposition}), we have
\begin{equation} \label{fx}
\te f(x)=\sum_{j=1}^p (\alpha_j \circ \phi) f(X_j).
\end{equation}
To prove the theorem, it suffices to show that $\phi$ is a local diffeomorphism and that
(\ref{dphi}) holds for
$x=X_j$, $j=1, \cdots, p$, and for any $ w \in D_{w_0}$.

Since $[X_i, X_j]=0$ on $O_{\phi(w_0)}$, where $i, j=1, \cdots, p$, we have
$[f(X_i), f(X_j)]$ $=0$ on $D_{w_0}$.
Thus there exist local coordinates $(w^1, \cdots, w^p)$ on an open neighborhood $\tilde{D}_{w_0} \subset D_{w_0}$ of
$w_0$ such  that $f(X_j)=\partial/\partial w^j$ on $\tilde{D}_{w_0}$ for $j=1, \cdots, p$.
On the open subset $O_{\phi(w_0)} \subset V$, we have $\left[X_j, v^i X_k \right]=\delta_{ji} X_k$, where $i, j, k=1, \cdots, p$,
and $\delta_{ji}$ is the Kronecker Delta.
By (\ref{fx}),
$$(f(X_j)(v^i \circ \phi))  f(X_k)=f\left([X_j, v^i X_k] \right)=\delta_{ji}f(X_k)$$
on $D_{w_0}$.
Thus $\partial(v^i \circ \phi)/\partial w^j=\delta_{ji}$ at $w_0$,
which implies that the differential of $\phi$ at $w_0$ is non-degenerate.
By appropriately shrinking $D_{w_0}$, we can assume that $\phi$ is invertible on $D_{w_0}$.
Let $\beta \in C^{\infty}(W)$. It follows from (\ref{fx}) that
\begin{eqnarray*}
   && (f(X_j)\beta) f(X_j) = [f(X_j), \beta f(X_j)]=f([X_j,  (\beta\circ \phi^{-1}) X_j])\\
  && =f((X_j(\beta\circ \phi^{-1})) X_j )=((X_j(\beta\circ \phi^{-1}))\circ \phi) f(X_j)
\end{eqnarray*}
on $D_{w_0}$. Hence $ f(X_j)\beta=(X_j(\beta\circ \phi^{-1}))\circ \phi =d\phi^{-1}(X_j) \beta$ on $D_{w_0}$.
Thus (\ref{dphi}) holds for
$x=X_j$, $j=1, \cdots, p$, and for each $w \in D_{w_0}$.
\end{proof}


\section{Homomorphisms of algebras of sections of bundles \label{ah}}

Any positive dimensional simple unital $\br$-algebra
is isomorphic to one of the $n \times n$ matrix algebras $M_n(R)$,
where $n \in \bn$ and $R=\br, \bc$, or the quaternion algebra $\bh$
\cite[Corollary 2.69]{bresar}. Furthermore,
any positive dimensional simple unital $\bc$-algebra
is isomorphic to $M_n(\bc)$, where $n \in \bn$ \cite[Corollary 2.66 and Theorem 2.61]{bresar}.
The only $\bff$-algebra automorphism of $\bff$ is the identity map.
Thus any
$\bff$-algebra bundle with typical fiber $\bff$ is trivial.
There exist nontrivial $\bff$-algebra bundles with typical fibers that are
positive dimensional simple unital $\bff$-algebras.
We say that an $\bff$-algebra is central if its center
is $\bff$. The $\br$-algebras $M_n(\br)$, $M_n(\bh)$ and the $\bc$-algebras $M_n(\bc)$
are all central \cite[Lemma 1.15]{bresar}. By the Skolem-Noether Theorem, every algebra automorphism of
a positive dimensional central simple algebra is an inner automorphism \cite[Section 1.6]{bresar}.
Hence the group of $\bc$-algebra automorphisms of $M_n(\bc)$, where $n=2, 3, \cdots$, is
$PGL(n, \bc)$. Let $G$ be a finite dimensional connected Lie group.
There is a bijective map from the set of equivalence classes of continuous principal $G$-bundles over $S^n$, where $n=2, 3, \cdots$,
to the homotopy group $\pi_{n-1}(G)$ \cite[Theorem 4.4.3]{na}.
Note that $\pi_1(PGL(n, \bc))=\bz_n$, $n=2, 3, \cdots$. Given a topologically nontrivial $C^{\infty}$ principal $PGL(n, \bc)$-bundle
on $S^n$,
the associated algebra bundle,  arising from the natural action of $PGL(n, \bc)$ on $M_n(\bc)$, is also nontrivial.

\begin{thm} \label{finitesp}
Suppose that $A_1$ (respectively $A_2$) is a positive dimensional unital $\bff$-algebra,
where $\bff=\br$ or $\bc$,
$\pi_{\sa_1, V}: \sa_1 \to V$ (respectively $\pi_{\sa_2, W}: \sa_2 \to W$) is a $C^{\infty}$ $\bff$-algebra bundle
with the typical fiber $A_1$ (respectively $A_2$),
where $\dim_{\br} V=1, 2, \cdots$ and $\dim_{\br} W=0, 1, \cdots$.
Let $f: \fcsao \to \ftsat$ (respectively $f: \fsao \to \ftsat$) be an $\bff$-algebra homomorphism.
If $E_{w} \circ f|_{\fcsao} \not\equiv \zero$ for each $w \in W$,
then the nonempty set $\cs_{E_{w} \circ f}$ contains at most $\dim_{\bff} A_2$ points.
  If in addition $A_1$ is simple, then there exists an injective
  $\bff$-algebra homomorphism $A_1 \to A_2$. Furthermore,
  if we assume that $A_1$ is simple and $A_1 = A_2$,
  then
  there exist
  an $\tilde{\cf}$ map $\phi: W \to V$ and an $\tilde{\cf}$ isomorphism of $\bff$-algebra bundles $\gamma: \phi^{\ast}\sa_1 \to \sa_2$
  covering the identity map on $W$
  such that
  \begin{equation} \label{algebrahomo}
  f(x)=\gamma \circ (\phi^{\ast} x), \hspace{2mm} x \in \fcsao \hspace{2mm} (\text{respectively} \hspace{2mm} x \in \fsao).
  \end{equation}
\end{thm}
\begin{proof}
By Proposition \ref{sp}, $\cs_{E_{w} \circ f} \not=\emptyset$ for each $w \in W$.
For any $a \in A_2$, we define the linear operator $l_a: A_2 \to A_2$ by $l_a(b)=ab$, $b \in A_2$.
Let $w_0 \in W$ and $v_0 \in \cs_{E_{w_0} \circ f}$. Take $x_0 \in \fcsao$ such that
$x_0$ is $\one$ on a neighborhood of $v_0$.
We can find $\tx_0 \in \fcsao$ such that $\supp \tx_0 \subset \{v \in V: x_0(v)=\one\}$
and $E_{w_0} \circ f(\tx_0) \not=\zero$. Since $x_0 \tx_0=\tx_0$,
we have $E_{w_0} \circ f(x_0) \not=\zero$, and the subspace $l_{E_{w_0} \circ f(x_0)}((\sa_2)_{w_0}) \subset (\sa_2)_{w_0}$
is not $\{\zero\}$.

Let $\ax_0 \in \fcsao$ be such that $\ax_0$ is $\one$ on a neighborhood of $v_0$
and $\supp \ax_0 \subset \{v \in V: x_0(v)=\one\}$.
Since $f(x_0) f(\ax_0) =f(\ax_0)$, any element in the subspace $l_{E_{w_0} \circ f(\ax_0)}((\sa_2)_{w_0}) \not=\{\zero\}$
is a fixed point of the linear operator $l_{E_{w_0} \circ f(x_0)}$.
Take pairwise distinct points $v_1, \cdots, v_k \in \cs_{E_{w_0} \circ f}$
and $x_1, \cdots, x_k \in \fcsao$ such that $x_j$ is $\one$
on a neighborhood of $v_j$ and $\supp x_i \cap \supp x_j =\emptyset$ for
$i \not=j$, where $i, j=1, \cdots, k$. Let $a_j \in (\sa_2)_{w_0} \setminus \{\zero\}$
be a fixed point of the map $l_{E_{w_0} \circ f(x_j)}$. Note that $x_i x_j \equiv\zero$ for $i \not=j$.
So
$a_j$ is in the kernel of $l_{E_{w_0} \circ f(x_i)}$. It is clear that $a_1, \cdots, a_k$ are linearly independent.
Thus $k \le \dim_{\bff} A_2$.

For any two-sided ideal $\ci \subset \fcsao$
and for any $v \in V$, $E_v(\ci)$ is a two-sided ideal of the fiber $(\sa_1)_v$.
Let  $\ci_{w_0}$ be the kernel
of $E_{w_0} \circ f|_{\fcsao}$.
For any $x_1 \in \fcsao$ with $x_1(v_0)=\one$, there exists an open
neighborhood $O$ of $v_0$ such that $x_1(v)$ is invertible
for each $v \in O$. Thus we can find another open neighborhood $O^{\ast} \subset O$ of $v_0$
and $x_2 \in \fcsao$ such that $x_1 x_2$ is $\one$ on $O^{\ast}$.
So $E_{w_0} \circ f(x_1 x_2) \not=\zero$
and $\one \not\in E_{v_0}(\ci_{w_0})$. If $A_1$ is simple, then $E_{v_0}(\ci_{w_0})=\{\zero\}$.
Take precompact open neighborhoods $O_{w_0}$, $O'_{w_0}$ of $\cs_{E_{w_0} \circ f}$ such that $\overline{O_{w_0}} \subset O'_{w_0}$
and there exists a local trivialization $\Psi_{w_0}: \pi^{-1}_{\sa_1, V}(O'_{w_0}) \to O'_{w_0} \times A_1$ of $\pi_{\sa_1, V}$.
For any $\alpha \in A_1$, choose $x_{\alpha} \in \fcsao$ such that $\supp x_{\alpha} \subset O'_{w_0}$ and
$\Psi_{w_0} \circ x_{\alpha}|_{O_{w_0}}$ is the constant section $\alpha$.
Define a map $\varrho_{w_0}: A_1 \to \fcsao$ by $\varrho_{w_0}(\alpha)=x_{\alpha}$. By Proposition \ref{sp},
$E_{w_0} \circ f \circ \varrho_{w_0}: A_1 \to (\sa_2)_{w_0} \simeq A_2$ is a well-defined
injective $\bff$-algebra homomorphism.

If $A_1$ is simple and $A_1 = A_2$, it follows from Theorem \ref{finiterankbundle} that
there exist
  an $\tilde{\cf}$ map $\phi: W \to V$ and an $\tilde{\cf}$ isomorphism of vector bundles $\gamma: \phi^{\ast}\sa_1 \to \sa_2$
  covering the identity map on $W$ such that (\ref{algebrahomo}) holds for any $x \in \fcsao$.
  Since $f$ is an $\bff$-algebra homomorphism, $\gamma$ is actually
  an isomorphism of $\bff$-algebra bundles.
  If the domain of $f$ is $\fsao$, then (\ref{algebrahomo}) holds for $x x^{\ast} \in \fcsao$,
 where $x \in \fsao$ and $x^{\ast} \in \fcsao$.
 This implies that (\ref{algebrahomo}) also holds for $x$.
\end{proof}

Let $\bff^k$ be the product of $k$ copies of the algebra $\bff$, where $k \in \bn$, and let
$\phi_1$, $\cdots$, $\phi_k$ $\in$ $C^{\infty}(W, V)$. Define $f_{k}$ as
the algebra homomorphism
$$C^{\infty}_c(V, \bff^k) \ni x=(x_1, \cdots, x_k) \mapsto (x_1 \circ \phi_1, \cdots, x_k \circ \phi_k) \in C^{\infty}(W, \bff^k).$$
It is clear that the upper bound in Theorem \ref{finitesp} for the number of elements in $\cs_{E_w \circ f_k}$ is sharp.
It follows from
Theorem \ref{finitesp} that $V$ is $\fcvf$-realcompact, even though it is not necessarily realcompact.
Furthermore,
a nonzero algebra homomorphism $f: \fvf \to \bff$
is the evaluation at some point of $V$ if and only if $f|_{\fcvf} \not\equiv 0$.

\begin{cor} \label{isoalgebra}
Suppose that $A_1$ (respectively $A_2$) is a positive dimensional simple unital $\bff$-algebra, where $\bff=\br$ or $\bc$,
and $\pi_{\sa_1, V}: \sa_1 \to V$ (respectively $\pi_{\sa_2, W}: \sa_2 \to W$) is a $C^{\infty}$ $\bff$-algebra bundle
with the typical fiber $A_1$ (respectively $A_2$).
 Let $f$ be an algebra isomorphism
  $\fcsao \to \fcsat$ (respectively $\fsao \to \fsat$). Then there exist
  a bi-$\cf$ map $\phi: W \to V$ and an $\cf$ isomorphism of algebra bundles $\gamma: \phi^{\ast}\sa_1 \to \sa_2$
  covering the identity map on $W$
  such that
  $f(x)=\gamma \circ (\phi^{\ast} x)$ for each $x \in \fcsao$ (respectively
  for each $x \in \fsao$).
\end{cor}
\begin{proof}
 The center of $A_i$ is either $\br$ or $\bc$, where $i=1, 2$.
We write $\pi_{\zsao, V}: \zsao \to V$ for the subbundle of $\pi_{\sa_1,V}$ such that the fiber of $\pi_{\zsao, V}$
   at $v \in V$ is the center of the fiber of $\pi_{\sa_1, V}$ at the same point.
   It is clear that $f(\cf_c(\pi_{\zsao, V})) = \cf_c(\pi_{Z(\sa_2), W})$
   (respectively $f(\cf(\pi_{\zsao, V}))$ $=$ $\cf(\pi_{Z(\sa_2), W})$). The bundles $\pi_{\zsao, V}$ and $\pi_{Z(\sa_2), W}$
   are trivial. If the domain of $f$ is $\fsao$ and
   $x \in \cf(\pi_{\zsao, V})$ is non-vanishing, then $x$ is invertible.
   Hence $f(x)$ is non-vanishing. By Theorem \ref{diffeo}, $E_w \circ f|_{\fcsao} \not\equiv \zero$
   for any $w \in W$. The conclusion of the corollary follows from Theorem \ref{finitesp}.
\end{proof}

By Corollary \ref{isoalgebra}, there exists an algebra isomorphism
$C^{\infty}_c(\pi_{\sa_1, V})$ $\to$ $C^{\infty}_c(\pi_{\sa_2, W})$ (respectively $C^{\infty}(\pi_{\sa_1, V})$ $\to$
$C^{\infty}(\pi_{\sa_2, W})$)
if and only if the algebra bundles $\sa_1$ and $\sa_2$ are isomorphic.


\end{document}